\documentclass[10pt]{article}
\usepackage{amsmath, amssymb,epsfig,bm}
\usepackage{latexsym}

\usepackage{graphics}
\usepackage{color}

\numberwithin{equation}{section}

\def\debproof{ {\bf Proof.} }
\def\finproof{\hfill {\small $\Box$} \\}

\newcommand{\bE}{{\bf E}}
\newcommand{\rE}{{\rm E}}
\newcommand{\bB}{{\bf B}}

\newcommand{\bH}{{\bf H}}
\newcommand{\bJ}{{\bf J}}
\newcommand{\bF}{{\bf F}}
\newcommand{\rF}{{\rm F}}

\newcommand{\bX}{{\bf X}}
\newcommand{\bn}{{\bf n}}
\newcommand{\bu}{{\bf u}}
\newcommand{\bv}{{\bf v}}
\newcommand{\bw}{{\bf w}}
\newcommand{\bx}{{\bf x}}

\newcommand{\curl}{{\bf curl}}
\newcommand{\re}{{\rm Re}}
\newcommand{\im}{{\rm Im}}
\newcommand{\bbf}{{\bf f}}
\newcommand{\be}{{\bf e}}
\newcommand{\Div}{{\rm Div}}

\newcommand{\q}{\quad}

\newtheorem{lem}{Lemma}[section]

\newtheorem{cor}{Corollary}[section]
\newtheorem{thm}{Theorem}[section]
\newtheorem{rem}{Remark}[section]

\newtheorem{alg}{Algorithm}[section]

\begin{document}

\title{Mathematical and numerical study of a three-dimensional \\ inverse eddy current problem 
}

\author{Junqing Chen\thanks{\footnotesize
Department of Mathematical Sciences, Tsinghua University, Beijing
100084, China. The work of this author was partially supported by China NSF under the grant 11871300, 91630205 and 11771440. (jqchen@tsinghua.edu.cn).}
\and Ying Liang\footnotemark[2] 
\and Jun Zou\thanks{\footnotesize Department of Mathematics, The Chinese University of Hong Kong, 
Shatin, Hong Kong. 
The work of J.\,Zou was substantially supported by Hong Kong RGC General Research Fund 
(Project 14304517) and National Natural Science Foundation of China/Hong Kong Research Grants Council Joint Research Scheme 2016/17 (Project N\_CUHK437/16). 
(yliang@math.cuhk.edu.hk, zou@math.cuhk.edu.hk).}
}

\date{}
\maketitle

\begin{abstract}
We study an inverse problem associated with an eddy current model.
We first address the ill-posedness of the inverse problem
by proving the compactness of the forward map with respect to the conductivity and 
the non-uniqueness of the recovery process. 
Then by virtue of non-radiating source conceptions, we establish a regularity result 
for the tangential trace of the true solution on the boundary, which is necessary to justify our subsequent 
mathematical formulation. 
After that, we formulate the inverse problem as a constrained optimization problem with an appropriate 
regularization and prove the existence and stability of the regularized minimizers. 
To facilitate the numerical solution of the nonlinear non-convex constrained optimization, 
we introduce a feasible Lagrangian and its discrete variant. Then the gradient of the objective functional is derived using the adjoint technique.  By means of the gradient, a nonlinear conjugate gradient method is formulated 
for solving the optimization system, and a Sobolev gradient is incorporated to accelerate the iterative process.
Numerical examples are provided to demonstrate the feasibility of the proposed algorithm.
\end{abstract}



\section{Introduction}\label{sect1}
Eddy current inversion is a challenging mathematical and numerical process, but it is 
one of the most popular nondestructive detection techniques. The inversion technique 
has attracted great attention in various important applications,
such as geophysical prospecting, flaw detection, safety inspection and biomedical imaging 
\cite{ABF,ACCGV,ACCVW,Haber,NB,Peyton,RCV,Zhdanov}. 
The eddy current method is based on the low frequency approximation of Maxwell's equation, and is much more sensitive to conductivity of materials 
when compared with the inversion by using full electromagnetic Maxwell system. 
There are two advantages by using the low frequency electromagnetic data in detection. First, low frequency electromagnetic wave can penetrate deeply in the lossy medium such as metal structure and the earth. It is well known that the intensity of electromagnetic wave will decay exponentially in lossy medium with respect to the penetration depth, and the intensity of higher frequency wave will decay faster \cite{Jackson}.
Second, the forward problem needs to be solved repeatedly in most inversion methods. 
While the full Maxwell's equations are difficult to solve numerically and efficiently themselves, 
the eddy current approximation of Maxwell's equations is a diffusion equation which can be solved 
with fast algorithms \cite{CCCZ,HX}. Therefore, the eddy current inversion method is widely used in nondestructive testing \cite{Peyton,RLL} and geophysical prospecting \cite{Haber, Zhdanov}.


Most inverse problems are known to be ill-posed.
We will study two important questions before we formulate our inverse model, i.e., 
the uniqueness and stability of the recovery. 
The analysis of these basic issues are very different with different inverse problems; 
see, e.g., \cite{FJZ} for time domain inverse Maxwell problem, \cite{CZ} for parameter identification problem with elliptic systems, and \cite{ABF, BN, RCV} for inverse Maxwell's source problems and inverse eddy current source problems. 
To the best of our knowledge, the uniqueness and stability analysis of the inverse eddy current problem have not 
been studied yet.
We shall investigate these two fundamental issues, then formulate and analyze the underlying constrained optimization problem as well as to propose some numerical method for 
the minimization. We start with the well-posedness of the forward eddy current problem, and 
establish a regularity result for the tangential trace of the true solution on the boundary by virtue of 
non-radiating source conceptions. This regularity is important to justify our usage 
of an appropriate selected misfit functional.
We then prove the compactness of the forward operator mapping the conductivity to the electric field 
and study the non-uniqueness of the inverse eddy current problem.  
With these preparations, we will formulate the ill-posed eddy current inverse problem 
into a nonlinear and non-convex constrained minimization with an appropriate regularization 
and show the existence and stability of the regularized minimizers. 
To facilitate the numerical solution of the nonlinear non-convex optimization constrained 
with the complex-valued eddy current model, 
we introduce a feasible Lagrangian and its discrete variant 
in terms of both real and imaginary parts of the constrained PDE.
Then we derive the gradient of the objective functional with the adjoint technique. 
For solving the nonlinear PDE constrained optimization, we formulate a nonlinear conjugate gradient (NLCG) method,  with the step size for the descent direction 
computed by a quadratic approximation to the state field. 
As the usual NLCG method converges very slowly, we incorporate a Sobolev preconditioner to improve 
the NLCG iteration.


The outline of the paper is as follows. In Section \ref{sect2}, we introduce the forward eddy current problem, 
present the well-posedness of the forward problem and prove the regularity of the the tangential trace of 
the true solution. In Section \ref{sect3}, an inverse problem with a well-defined misfit functional is formulated 
and the ill-posedness of the inversion problem is investigated. 
Then we add a regularization term to the optimization objective functional and prove the existence and 
stability of the minimizers. In Section \ref{sect4}, we first introduce a Lagrangian
associated with the regularized optimization problem, then introduce the gradient of the objective functional with adjoint technique, and further study the properties of the adjoint state equation.  Moreover, 
the finite element discretization of the optimization problem is also formulated and studied 
in the same section, and a nonlinear conjugate gradient method is proposed for 
the optimization system.  
We show some numerical examples in Section \ref{sect5} to illustrate the feasibility of the proposed algorithm, and 
present some concluding remarks in Section \ref{sect6}.

\section{The forward problem}\label{sect2}
In this section, we introduce the forward model for eddy current inversion and present some necessary preliminaries.
The eddy current equation is the low frequency approximation of Maxwell's equation. 
As we mentioned in the Introduction, the eddy current field can penetrate deeply in conducting materials.
Moreover, as an electromagnetic method, this method can distinguish the conductor (metal, water) from the insulator (oil, rock) and is an important modality in nondestructive detection. 
The eddy current problem has been studied extensively in the literature \cite{RV}. 
The governing equations for the forward problem read
\begin{eqnarray*}
\left\{
\begin{array}{rlll}
\nabla\times \bE &=& i \omega\mu \bH &\mbox{ in }\mathbb{R}^3,\\
\nabla\times \bH &=& \sigma^\prime \bE + \bJ_s &\mbox{ in } \mathbb{R}^3, 
\end{array}
\right.
\end{eqnarray*}
where $\bE$, $\bH$ are electric and magnetic fields respectively, $\mu$ is the magnetic permeability, $\sigma^\prime$ is the conductivity of the medium and $\bJ_s$ is the source current. 
While the equations hold in the whole space $\mathbb{R}^3$, we consider the problem in a bounded domain $\Omega\subset \mathbb{R}^3$ as in many applications and theories, 
and boundary conditions are specified later to form a well-posed problem.

\begin{figure}
\centering{
\includegraphics[width=0.5\textwidth]{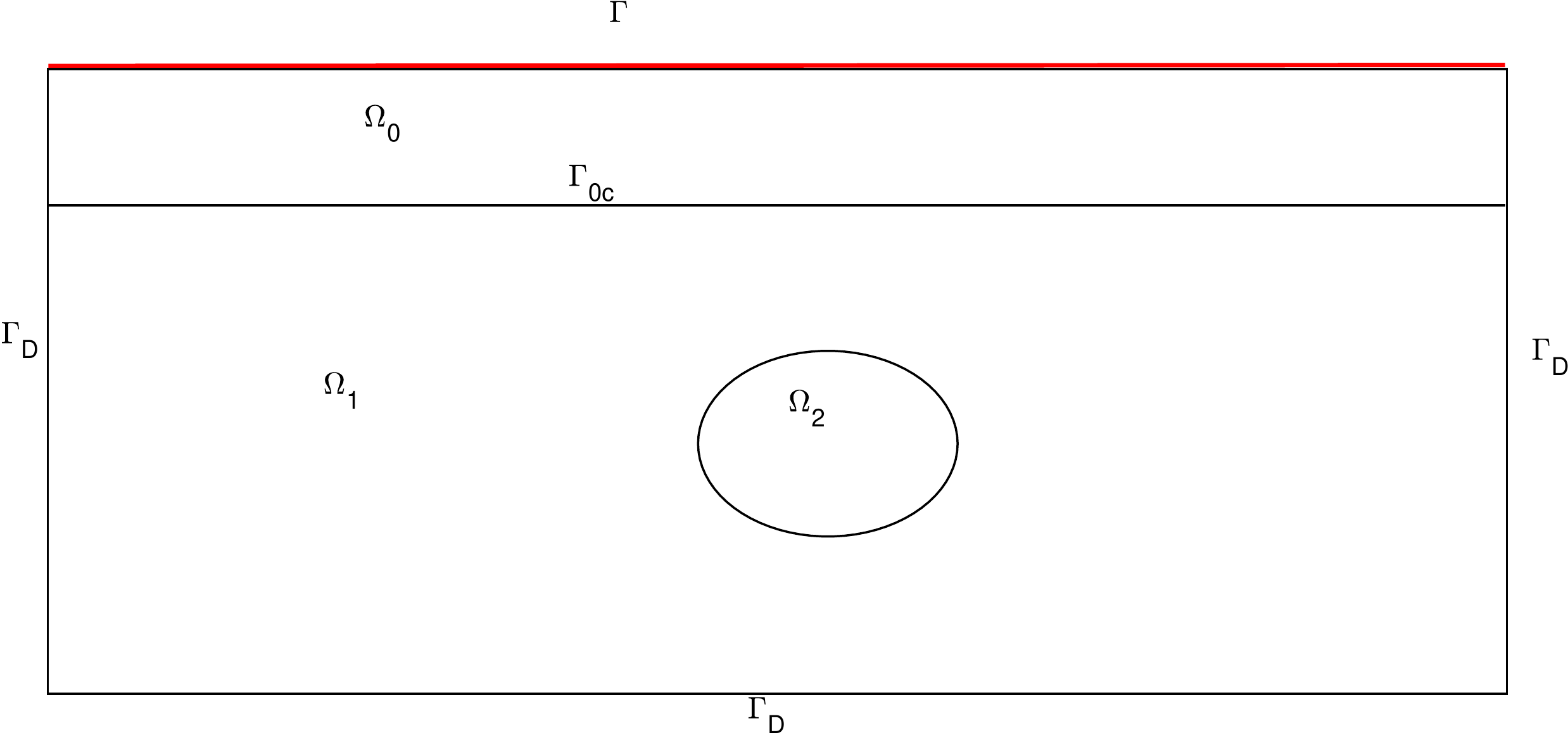}
}
\caption{The geometric setting of the problem}\label{conf1}
\end{figure}

Now we start with some assumptions for the further consideration of the eddy current model. 
In the rest of this paper, we concentrate on the electric acquisition case, that is, the measurement data is collected 
for the tangential components of the electric field on $\Gamma$, part of the boundary $\partial\Omega$. 
We assume that $\Omega$ is a convex domain, with a piecewise smooth boundary and 
a simply-connected subdomain $\Omega_0$ occupied by air, 
hence the conductivity $\sigma^\prime$ vanishing  in $\Omega_0$. 
A typical geometric setting of the problem in a 2D cross-section is shown in Figure \ref{conf1}, where $\overline\Omega=\overline\Omega_0\cup\overline\Omega_1\cup\overline\Omega_2$. 
The material parameter is a different function at each subdomain. We write 
the conductivity $\sigma^\prime(x)$ in $\Omega$ as 
\begin{eqnarray*}
\sigma^\prime(x)=\sigma_0+\sigma(x), 
\end{eqnarray*}
where $\sigma_0$ is the constant background conductivity which is supported in $\Omega_1\cup\overline\Omega_2$ and known a priori. $\sigma(x)$ is the abnormal conductivity. Both $\sigma(x)$ and its support 
$\Omega_2$ are unknown and are our target to recover simultaneously. 
We shall write $\Omega\setminus\overline\Omega_0$ as $\Omega_c$, 
then $\sigma_0+\sigma(x)$ is supported in $\Omega_c$. We further assume that $\sigma(x)$ 
is compactly supported in $\Omega_c$. The interface between $\Omega_0$ and $\Omega_c$ is denoted by $\Gamma_{0c}$ and assumed to be a simply-connected Lipschitz polyhedral surface, with both 
domains $\Omega_0$ and $\Omega_c$ being polyhedrons and simply-connected. 
In our subsequent study, $\mu$ is assumed to be piecewise constant physically, 
and the source $\bJ_s$ is compactly supported in $\Omega_0$, and  $\nabla\cdot\bJ_s=0$.

\subsection{The $\bE$-based eddy current model and its inverse problem}
 By eliminating $\bH$ in the eddy current equations, we derive the electric field system
\begin{eqnarray}
\left\{\begin{array}{rlll}
\nabla\times(\mu^{-1} \nabla\times \bE)-i\omega(\sigma_0+\sigma) \bE&=&i\omega\bJ_s &\mbox{ in } \Omega\,, \\
\nabla\cdot\varepsilon \bE&=&0 &\mbox{ in }\Omega_0\,,
\end{array}\right.  \label{MTE}
\end{eqnarray}
which are complemented by with the interface condition
\begin{eqnarray}
[\mu^{-1}\bn\times\nabla\times\bE]=0 \mbox{ on } \Gamma_{0c}\cup \partial\Omega_2\, ,\label{MTE_ic}
\end{eqnarray}
and the boundary conditions
\begin{eqnarray}
\bn\times \nabla\times \bE=0 ~~\mbox{on} ~~\Gamma\,; \q 
\bn\cdot\bE=0 ~~\mbox{on} ~~\Gamma\,; \q 
\bn\times \bE=0 ~~\mbox{on} ~~\Gamma_D=\partial\Omega\setminus\overline\Gamma\,.
\label{BCs}
\end{eqnarray}
Here and in the sequel, $\bn$ denotes the outward normal on $\partial\Omega$.
We add a divergence free equation in the system \eqref{MTE} to ensure the uniqueness of the solution since  $\sigma^\prime=0$ in $\Omega_0$. The piecewise constant $\varepsilon$ is the electric permittivity in $\Omega_0$. The divergence free condition makes the field $\bE$ to be electric field in domain $\Omega_0$.
The surface $\Gamma$ is where we measure the data, i.e., the tangential components $\bn\times \bE$ 
of the electric field. The inverse eddy current problem of our interest is formulated as follows: 

\smallskip
%
Given the observation data $\bn\times \bE^{obs}$ on the measurement surface $\Gamma$, recover 
the conductivity distribution $\sigma(x)$ and its support $\Omega_2$.


 
\subsection{The weak formulation and regularity of the solution}
For the variational formulation of the electric field problem \eqref{MTE} and its well-posedness, 
we introduce the Sobolev spaces:
\begin{eqnarray*}
H_{\Gamma}(\curl;\Omega) & =& \big\{\bu\in L^2(\Omega)^3 ~\big|~\nabla\times \bu\in L^2(\Omega)^3, ~\bn\times \bu=0 \mbox{ on }\Gamma_D\big\},\\
H^1_{\Gamma}(\Omega_0)&=&\big\{v\in L^2(\Omega_0)~\big|~\nabla v\in L^2(\Omega_0)^3, ~v|_{\partial\Omega_0\setminus\Gamma}=0\big\},\\
\mathbf{Y}&=&\big\{\bu\in H_{\Gamma}(\curl;\Omega)~\big|~(\varepsilon \bu,\nabla\phi)=0 ~~\forall\,\phi\in H^1_\Gamma(\Omega_0)\big\},
\end{eqnarray*}
and the tangential trace space of $H_\Gamma(\curl;\Omega)$ on $\Gamma$:
 $$H^{-1/2}(\Div;\Gamma)=\big\{\bbf\in H^{-1/2}(\Gamma)^3 ~\big|~ \exists \ \bu \in H_\Gamma(\curl;\Omega) \  
 \mbox{such that} \ \bn\times\bu =\bbf  \big\},$$
 or equivalently \cite{Monk},
 $$H^{-1/2}(\Div;\Gamma)=\big\{\bbf\in H^{-1/2}(\Gamma)^3 ~\big|~ 
 \bn\cdot\bbf=0\ a.e.\  {\rm on}\ \Gamma;~ \Div_\tau \bbf\in H^{-1/2}(\Gamma) \big\}.$$
 Here $\Div_\tau $ is the surface divergence operator which will be formally defined on smooth surface  in Section \ref{sect4}.
We define a sesquilinear form $a: H_{\Gamma}(\curl; \Omega)\times H_{\Gamma}(\curl; \Omega)\rightarrow\mathbb{C}$ as
\begin{equation}\label{eq:a(u,v)}
a(\bE,\bF)=\int_\Omega \mu^{-1}\nabla\times \bE\cdot\nabla\times \overline\bF  -i\omega(\sigma+\sigma_0) \bE\cdot\overline\bF dx  \quad \forall\, \bE,\bF\in H_{\Gamma}(\curl; \Omega),
\end{equation}
where $\overline\bF$ denotes the vector-valued complex conjugate of $\bF$.
Then the weak formulation of problem \eqref{MTE} is:
Find $\bE\in\mathbf{Y}$ such that
\begin{eqnarray}
a(\bE,\bF)=i\omega\int_{\Omega}\bJ_s\cdot \overline\bF dx ~~\forall\, \bF\in\mathbf{Y}\label{MTE0}.
\end{eqnarray}
The following lemma implies that the well-posedness of the problem \eqref{MTE0}.
\begin{lem}
The problem \eqref{MTE0} has a unique solution $\bE\in\mathbf{Y}$.
\end{lem}
\debproof
The uniqueness is due to the fact that sesquilinear form $a(\cdot,\cdot)$ is coercive in space $\mathbf{Y}$. The proof of coercivity is similar to \cite{CCCZ}. For completeness, we sketch a proof here.
For any $\bE\in\mathbf{Y}$, $\bn\times \bE|_{\Gamma_{0c}}\in H^{-1/2}(\Div; \Gamma_{0c})$ \cite{BCS}. Let $H_\Gamma(\curl; \Omega_0)=\{\bu\in H(\curl; \Omega_0)~\big|~\bn\times \bu=0\mbox{ on } \partial\Omega_0\setminus\Gamma\}$.  By  Lax-Milgram theorem, there exists a unique  $\bB \in H(\curl; \Omega_0)$, $\bn\times \bB=0$ on $\Gamma_D$, $\bn\cdot \bB=0$ on $\Gamma$ and $\bn\times \bB=\bn\times \bE$ on $\Gamma_{0c}$, such that
\begin{eqnarray*}
\int_{\Omega_0}\nabla\times \bB\cdot\nabla\times \overline\bu dx+\int_{\Omega_0}\varepsilon \bB\cdot \overline \bu dx = 0 ~~\forall\, \bu\in H_\Gamma(\curl; \Omega_0).
\end{eqnarray*}
By the trace theorem,
\begin{eqnarray*}
\|\bB\|_{H(\curl;\Omega_0)}\leq C\|\bn\times \bB\|_{H^{-1/2}(\Div;\Gamma_{0c})}
=C\|\bn\times \bE\|_{H^{-1/2}(\Div;\Gamma_{0c})}
\leq C\|\bE\|_{H(\curl;\Omega_c)}.
\end{eqnarray*}
Moreover, we have
$$\nabla\cdot \varepsilon \bB =0 \mbox{ and } \nabla\cdot\varepsilon(\bE-\bB)=0 \mbox{ in }\Omega_0, $$
and $\bn\times(\bE-\bB)=0$ on $\partial\Omega_0\setminus\Gamma$ and $\bn\cdot (\bE-\bB)=0$ on $\Gamma$. 
Then we know 
\begin{eqnarray*}
\|\bE-\bB\|_{L^2(\Omega_0)}\leq C\|\nabla\times(\bE-\bB)\|_{L^2(\Omega_0)}, 
\end{eqnarray*}
and furthermore, 
\begin{eqnarray*}
\|\bE\|_{L^2(\Omega_0)}&\leq& C(\|\bE-\bB\|_{L^2(\Omega_0)}+\|\bB\|_{L^2(\Omega_0)})\\
&\leq& C(\|\nabla\times \bE\|_{L^2(\Omega_0)}+\|\bE\|_{H(\curl;\Omega_c)})\\
&\leq& C(\|\nabla\times \bE\|_{L^2(\Omega)}+\|\bE\|_{L^2(\Omega_c)}).
\end{eqnarray*}
This implies that
$|a(\bE,\bE)|\geq C \|\bE\|^2_{H(\curl;\Omega)}$
for all $\bE\in\mathbf{Y}$.
\finproof

It is difficult to solve problem \eqref{MTE0} numerically since it is hard to construct a conforming finite element  space of $\mathbf{Y}$. Therefore we reformulate the weak formulation \eqref{MTE0} as a saddle point problem by introducing a Lagrange multiplier to deal with the divergence-free condition in domain $\Omega_0$. 
The saddle point formulation of equation \eqref{MTE} reads: 
Find $(\bE,\phi)\in H_\Gamma(\curl;\Omega)\times H^1_\Gamma(\Omega_0)$ such that
\begin{equation}
\left\{
\begin{array}{rlll}
a(\bE,\bF)+b(\nabla\phi,\bF) &=&i\omega\int_{\Omega} \bJ_s\cdot\overline\bF dx& ~~\forall\,\bF\in H_\Gamma(\curl;\Omega), \\
b(\bE,\nabla\psi)&=&0& ~~\forall\, \psi\in H^1_\Gamma(\Omega_0).
\end{array}
\right.
\label{MTE-w}
\end{equation}
where $b: H_{\Gamma}(\curl;\Omega)\times H_{\Gamma}(\curl;\Omega)\rightarrow\mathbb{C}$ 
is a sesquilinear form given by 
$$
b(\bE,\bF)=\int_{\Omega_0}\varepsilon\bE\cdot\overline\bF dx~~\forall \,\bE,\bF\in H_{\Gamma}(\curl;\Omega).
$$

\begin{lem}[Uniqueness]\label{uniq}
There is at most one solution to \eqref{MTE-w}.
\end{lem}
\debproof
We only need to show that $\bE=0$ in $\Omega$ and $\phi=0$ in $\Omega_0$ provided $\bJ_s=0$.
First, we take $\bF$ as the zero extension of $\nabla\phi$ from $\Omega_0$ to $\Omega$, that is to say,
$$\bF=0 \q \mbox{ in}\q \Omega_c \q \mbox{and} \q 
\bF=\nabla\phi \q \mbox{in}\q \Omega_0,$$
which implies $\bF\in H_\Gamma(\curl;\Omega)$. 
We plug $\bF$ into the first equation of \eqref{MTE-w}, along with $\bJ_s=0$,  to get 
\begin{eqnarray*}
\int_{\Omega_0}|\nabla\phi|^2dx=0.
\end{eqnarray*}
So $\nabla\phi=0$ in $\Omega_0$, and by the boundary condition on $\partial\Omega_0\setminus\Gamma$, we have $\phi=0$.
Second, taking $\bF=\bE$, $\psi=\phi$ in \eqref{MTE-w}, we obtain 
\begin{eqnarray*}
a(\bE,\bE)=\int_{\Omega}\mu^{-1}|\nabla\times\bE|^2dx-i\omega\int_{\Omega_c}(\sigma+\sigma_0)|\bE|^2dx=0.
\end{eqnarray*}
This implies $\nabla\times\bE=0$ in $\Omega$ and $\bE=0$ in $\Omega_c$. By the tangential continuity of $\bE$, we know that
\begin{eqnarray*}
\left\{
\begin{array}{rlll}
\nabla\times \bE&=&0 &\mbox{ in } \Omega_0,\\
\nabla\cdot\bE&=&0& \mbox{ in } \Omega_0,\\
\bn\times\bE&=&0& \mbox{ on }\partial\Omega_0\setminus\Gamma,\\
\bn\cdot\bE&=&0& \mbox{ on } \Gamma.
\end{array}
\right.
\end{eqnarray*}
By the assumption, $\varepsilon$ is constant in the simply-connected domain $\Omega_0$. Then there exists $p\in H^1(\Omega_0)$ such that $\bE=\nabla p$ and
\begin{eqnarray*}
\left\{
\begin{array}{rlll}
\Delta p&=& 0 &\mbox{ in } \Omega_0\\
\frac{\partial p}{\partial\bn}&=&0& \mbox{ on } \Gamma\\
p&=&C& \mbox{ on } \partial\Omega_0\setminus\Gamma
\end{array}
\right.
\end{eqnarray*}
for some constant $C$.
It is easy to know that the unique solution of the above system is $p=C$, so $\bE=0$ in $\Omega_0$. 
This completes our proof.\finproof

\begin{thm} \label{EXIST}
The equation \eqref{MTE-w} has a unique solution $(\bE,\phi)\in H_\Gamma(\curl,\Omega)\times H^1_\Gamma(\Omega_c)$ and $\bE$ satisfies \eqref{MTE0}. Moreover, the following stability estimate holds:
\begin{eqnarray}
\|\bE\|_{H(\curl;\Omega)}+\|\phi\|_{H^1(\Omega_0)}\leq C\|\bJ_s\|_{L^2(\Omega)},\label{regular1}
\end{eqnarray}
where $C$ is a constant independent of $\bE$ and $\phi$.
\end{thm}
\noindent\debproof
The existence can be established by proving the equivalence between \eqref{MTE0} and \eqref{MTE-w}.
Let $\bE$ be the solution of \eqref{MTE0}, then it is clear that $\bE$ satisfies the second equation of \eqref{MTE-w}. If we can prove that there exists $\phi\in H^1_\Gamma(\Omega_0)$ such that $\bE$ and $\phi$ satisfy the first equation of \eqref{MTE-w},  by the uniqueness of solution to \eqref{MTE-w} we may conclude the existence of solution of \eqref{MTE-w}.

Now, for any $\bF\in H_\Gamma(\curl,\Omega)$, we can find a $\psi$ that satisfies
\begin{eqnarray*}
\int_{\Omega_0}\nabla\psi\cdot \nabla\overline\xi dx=\int_{\Omega_0}\bF\cdot\nabla\overline\xi dx ~~\forall\,\xi\in H^1_\Gamma(\Omega_0),
\end{eqnarray*}
then $\nabla\cdot(\bF-\nabla\psi)=0$ in $\Omega_0$. Let $\tilde\psi$ be an extension of $\psi$, 
 $$\tilde\psi=\left\{
 \begin{array}{cc}\psi &\mbox{in }\Omega_0,\\ 0 & \mbox{ otherwise. }\end{array}\right.$$
Let $A=\bF-\nabla\tilde\psi\in\mathbf{Y}$. Since $\tilde{\psi}$ is supported in $\Omega_0$, we have
\begin{eqnarray*}
a(\bE,\bF)+b(\nabla\phi,\bF)&=&a(\bE,A+\nabla\tilde\psi)+b(\nabla\phi,A+\nabla\tilde\psi)\\
&=&a(\bE,A)+a(\bE,\nabla\tilde\psi)+b(\nabla\phi,A)+b(\nabla\phi,\nabla\psi)\\
&=&a(\bE,A)+b(\nabla\phi,\nabla\psi)\\
&=&i\omega\int_{\Omega}\bJ_s\cdot \overline Adx + b(\nabla\phi,\nabla\psi).
\end{eqnarray*}
 The right-hand side of the first equation of \eqref{MTE-w} becomes
\begin{eqnarray*}
i\omega\int_{\Omega}\bJ_s\cdot \overline\bF dx=i\omega\int_{\Omega}\bJ_s \cdot \overline Adx+i\omega\int_{\Omega_0}\bJ_s \cdot\nabla\overline\psi dx.
\end{eqnarray*}
Let $\phi \in H^1_\Gamma(\Omega_0)$ be a solution to the variational system 
\begin{eqnarray}
b(\nabla\phi,\nabla\psi)= i\omega\int_{\Omega_0}\bJ_s \cdot\nabla\overline \psi dx ~~\forall\,\psi \in H^1_\Gamma(\Omega_0) \label{pphi}.
\end{eqnarray}
We know there exists a unique $\phi$ satisfying \eqref{pphi}. Actually $\phi=0$ because 
$\bJ_s$ is divergence-free and compactly supported in $\Omega_0$. 
With $\bE,\phi$ satisfying \eqref{MTE0} and \eqref{pphi} respectively, we have
\begin{eqnarray*}
a(\bE,\bF)+b(\nabla\phi,\bF)=i\omega\int_{\Omega}\bJ_s\cdot\overline\bF dx\, \quad 
\forall\,\bF\in H_\Gamma(\curl,\Omega)\,.
\end{eqnarray*}
Then $(\bE,\phi)$ is a solution to \eqref{MTE-w}. We can now conclude the existence and uniqueness of solution to \eqref{MTE-w} by Lemma \ref{uniq}. Furthermore, if $(\bE,\phi)$ is a solution to \eqref{MTE-w}, 
we readily see $\bE$ is a solution to \eqref{MTE0}.
\finproof

It is known that the tangential trace space of $H_\Gamma(\curl,\Omega)$ is
$H^{-1/2}(\Div;\Gamma)$ \cite{BCS}, 
i.e., $\bn\times\bE|_{\Gamma}\in H^{-1/2}(\Div;\Gamma)$ for all $\bE\in H_\Gamma(\curl;\Omega)$.  Let $\bn\times\bE^{obs}$ be the data on $\Gamma$, 
and $\bn\times\bE$ be the corresponding tangential part of the electric field $\bE$ 
on $\Gamma$ associated with the conductivity $\sigma$.
 Then a direct choice of the misfit of prediction is $\|\bn\times(\bE^{obs}-\bE)\|_{H^{-1/2}(\Div;\Gamma)}$. Unfortunately, this trace space is naturally equipped with the norm
$$\|\mathbf{f}\|_{H^{-1/2}(\Div;\Gamma)}=\inf_{\bu\in H_\Gamma(\curl,\Omega),\bn\times \bu=\mathbf{f}\mbox{ on }\Gamma}\|\bu\|_{H(\curl,\Omega)},$$
which is difficult to realize numerically.
It would be very convenient and important numerically if a computable norm, such as the $L^2$-norm on $\Gamma$, 
can be used for the recovery process.
Next, we demonstrate that the true solution $\bE$ to the problem \eqref{MTE} 
indeed have a higher regularity, suggesting us a computable norm on $\Gamma$.
\begin{thm}\label{regularity}
Assuming that $\Omega_0$ and $\Omega_c$ are polyhedral domains and $\Omega$ is convex, $\sigma_0$ is a constant in $\Omega_c$,  $\bE^{obs}$ is the solution to \eqref{MTE} with the exact conductivity $\sigma_0+\sigma_e$, then for any $\sigma$ we have
\begin{eqnarray}
(\bE(\sigma)-\bE^{obs})|_{\Omega_0}\in H^{1/2}(\Omega_0).
\end{eqnarray}
\end{thm}
\debproof
It follows from \eqref{MTE} that 
\begin{eqnarray*}
\nabla\times(\mu^{-1}\nabla\times\bE(\sigma))-i\omega(\sigma_0+\sigma)\bE(\sigma)=i\omega\bJ_s,\\
\nabla\times(\mu^{-1}\nabla\times\bE^{obs})-i\omega(\sigma_0+\sigma_e)\bE^{obs}=i\omega\bJ_s,
\end{eqnarray*}
from which we can easily deduce
\begin{eqnarray}
\nabla\times(\mu^{-1}\nabla\times(\bE(\sigma)-\bE^{obs}))-i\omega\sigma_0(\bE(\sigma)-\bE^{obs})
=\bJ_e,
\end{eqnarray}
where $\bJ_e=i\omega\sigma\bE(\sigma)-i\omega\sigma_e\bE^{obs}$.
For $\bJ_e$, we have the following decomposition, 
\begin{eqnarray*}
\bJ_e=\bJ_0+\nabla\phi,
\end{eqnarray*}
where $\nabla\cdot\bJ_0=0$ and $\Delta\phi=\nabla\cdot\bJ_e$,  $\phi\in H^1_0(\Omega_c)$. Then we let $\bE(\sigma)-\bE^{obs}=\bE_r+\bE_\phi$ and $\bE_r$, $\bE_\phi$ satisfy the following two systems respectively, with the boundary conditions \eqref{BCs}, i.e.,

\begin{eqnarray*}
\left\{
\begin{array}{rlll}
\nabla\times(\mu^{-1}\nabla\times\bE_r)-i\omega\sigma_0\bE_r&=&\bJ_0 & \mbox{ in }\Omega,\\
\nabla\cdot\varepsilon\bE_r&=&0 & \mbox{ in }\Omega_0,
\end{array}
\right.
\end{eqnarray*}
and
\begin{eqnarray*}
\left\{
\begin{array}{rlll}
\nabla\times(\mu^{-1}\nabla\times\bE_\phi)-i\omega\sigma_0\bE_\phi&=&\nabla\phi & \mbox{ in }\Omega,\\
\nabla\cdot\varepsilon\bE_\phi&=&0 & \mbox{ in }\Omega_0.
\end{array}
\right.
\end{eqnarray*}
By the assumption on $\Omega_0$ and $\Omega_c$, we know that $\Gamma_{0c}$ is a Lipschitz interface. With the help of Theorem 6.1 in \cite{CDN1}, we find that $\bE_r|_{\Omega_0} \in H^{1/2}(\Omega_0)$. As for $\bE_\phi$, with the arguments in Theorem \ref{NRsource} of Section \ref{sect3}, we know that
$\nabla\phi$ is a non-radiating source, then $\bE_\phi|_{\Omega_0}=0$.  Then we complete the proof by noting that
$ 
\bE(\sigma)-\bE^{obs}=\bE_r
$ 
on $\Omega_0$.
\finproof

Theorem\,\ref{regularity} 
implies that the regularity of the solution to \eqref{MTE0} in subdomain $\Omega_0$ is higher than the global regularity. With this result, we further derive the following estimate.
\begin{lem}\label{regularity1}
With the same assumptions and notations as in Theorem \ref{regularity}, we have the estimate
\begin{eqnarray*}
\|\bE(\sigma)-\bE^{obs}\|_{H^{1/2}(\Omega_0)}\leq C\|\bJ_e\|_{L^2(\Omega_c)},
\end{eqnarray*}
where $C$ is independent of $\sigma$.
\end{lem}
\debproof
Recall the definition of $\bE_r$ and $\bE_\phi$ in Theorem \ref{regularity}. Let us denote
$$X_N(\Omega_0)=\big\{\bu,\nabla\times\bu\in L^2(\Omega_0)^3, \nabla\cdot\varepsilon\bu\in L^2(\Omega_0),\bn\times\bu=0 \mbox{ on }\Gamma_D, \bn\cdot\bu=0\mbox{ on }\Gamma\big\}.$$
Then, $\bE_r|_{\Omega_0}\in X_N(\Omega_0)$. From the imbedding of $X_N(\Omega_0)$ in $H^{1/2}(\Omega_0)$ as in Theorem 6.1 of \cite{CDN1}, we know that
\begin{eqnarray*}
\|\bE_r\|_{H^{1/2}(\Omega_0)}\leq C(\|\bE_r\|^2_{L^2(\Omega_0)^3}+\|\nabla\times\bE_r\|^2_{L^2(\Omega_0)^3}+\|\nabla\cdot\varepsilon\bE_r\|^2_{L^2(\Omega_0)})^{1/2}.
\end{eqnarray*}
Since $\bE_\phi=0$ and $\nabla\cdot\varepsilon\bE_r=0$ in $\Omega_0$, we have with the help of the estimate \eqref{regular1} that
\begin{eqnarray*}
\|\bE(\sigma)-\bE^{obs}\|_{H^{1/2}(\Omega_0)}\leq C(\|\bE(\sigma)-\bE^{obs}\|^2_{L^2(\Omega_0)^3}+\|\nabla\times(\bE(\sigma)-\bE^{obs})\|^2_{L^2(\Omega_0)^3})^{1/2}
\leq C\|\bJ_e\|_{L^2(\Omega_c)^3}.
\end{eqnarray*}
\finproof
\begin{rem}
As pointed out in \cite{ABDG}, provided that the domain $\Omega_0$ is of class $C^{1,1}$ or convex, $X_N\cap H_0(\curl;\Omega_0)$ is continuously imbedded in $H^1(\Omega_0)^3$. In our present situation, $\Omega_0$ is a Lipschitz polyhedron domain, and $X_N(\Omega_0)$ is imbedded in $H^{s}(\Omega_0)^3$ for some $s>1/2$. These imbedding results play an important role in the subsequent sections of the paper.
\end{rem}
%
\section{Ill-posedness of the inverse problem}\label{sect3}
In this section, we investigate the ill-posedness of the eddy current inverse problem. 
We know the solution $(\bE,\phi)\in H_\Gamma(\curl;\Omega)\times H^1_\Gamma(\Omega_0)$ 
to problem \eqref{MTE-w} depends on the conductivity 
$\sigma_0+\sigma(x)$. But in the setting of our inverse problem, $\sigma_0$ is known, so we shall 
write $\bE(\sigma)$ to emphasize its dependence on $\sigma$. 
The ill-posedness of the eddy current inverse problem is basically determined by 
the nature of the forward operator $\bE(\sigma)$. 

%
\subsection{Compactness of the forward operator $\bn\times\bE(\sigma)$}
We first present a result about the continuity of $\bn\times\bE(\sigma)$.
\begin{lem}\label{lem:3.1}
For any sequence $\{\sigma_n\}\subset H^1_0(\Omega_c)$ such that 
$\sigma_n\rightarrow \sigma_*$ in $L^2(\Gamma)$ as $n\rightarrow\infty$, it holds that 
\begin{eqnarray}
\lim_{n\rightarrow\infty}\|\bn\times(\bE(\sigma_n)-\bE(\sigma_*))\|_{L^2(\Gamma)}=0.
\end{eqnarray}
\end{lem}
\debproof
By the definition, $(\bE(\sigma_n),\phi(\sigma_n))$ and $(\bE(\sigma_*),\phi(\sigma_*))$ are the solutions to equation \eqref{MTE-w} with $\sigma$ replaced by $\sigma_n$ and $\sigma_*$, respectively. Then let $\hat\bE_n=\bE(\sigma_n)-\bE(\sigma_*)$ and $\hat\phi_n=\phi(\sigma_n)-\phi(\sigma_*)$, it is easy to check that  $(\hat\bE_n,\hat\phi_n)$ satisfies \eqref{MTE-w} with $\sigma_n$ and $(\sigma_n-\sigma_*)\bE(\sigma^*)$ in places of $\sigma$ and $\bJ_s$, respectively. 

Then from Theorem \ref{EXIST} we can deduce that
\begin{eqnarray*}
\|\hat\bE_n\|_{H(\curl;\Omega)}+\|\hat\phi_n\|_{H^1(\Omega_0)}&\leq& C\|(\sigma_n-\sigma_*)\bE(\sigma_*)\|_{L^2(\Omega)^3}\\
&\leq &C\|\sigma_n-\sigma_*\|_{L^2(\Omega_c)}\|\bE(\sigma_*)\|_{L^2(\Omega_c)^3}\\
&\leq &C\|\sigma_n-\sigma_*\|_{L^2(\Omega_c)}\|\bJ_s\|_{L^2(\Omega)^3}.
\end{eqnarray*}
With the same argument as in Lemma \ref{regularity1}, we know $\hat\bE_n|_{\Omega_0}\in H^{1/2}(\Omega_0)$, and hence
\begin{eqnarray*}
\|\hat\bE_n\|_{H^{1/2}(\Omega_0)}\leq C\|\sigma_n-\sigma_*\|_{L^2(\Omega_c)}\|\bJ_s\|_{L^2(\Omega)^3}.
\end{eqnarray*}
Then by  trace theorem we can conclude the convergence result.
\finproof

By Theorem \ref{EXIST}, given a parameter $\sigma\in H^1_0(\Omega_c)$, there exists a solution $\bE(\sigma)$ of \eqref{MTE-w}, and $\bE(\sigma)$ determines the tangential component 
$\bn\times\bE(\sigma)$ on $\Gamma$. In this way one can define the forward map $\bE(\sigma)$ from $H^1_0(\Omega_c)$ to $L^2(\Gamma)^3$.  In the following lemma, we prove that this implicit map from parameter $\sigma$ to the tangential field $\bn\times \bE(\sigma)$ is compact, 
and this implies that the eddy current inverse model 
that uses the data on $\Gamma$ to determine $\sigma$ is ill-posed.
\begin{lem}
The map from $\sigma\in H^1_0(\Omega_c)$ to $\bn\times\bE(\sigma)|_\Gamma\in L^2(\Gamma)$ 
is compact.
\end{lem}
\debproof
Let $\{\sigma_n\}_{n=1}^\infty$ be a bounded sequence in $ H^1_0(\Omega_c)$. Since $H^1(\Omega_c)$ is compactly imbedded in $L^2(\Omega_c)$, there is a subsequence, still denoted as $\{\sigma_n\}$, that converges to  $\sigma^*$ in $L^2(\Omega_c)$.
By Lemma \ref{regularity}, for each $\sigma_n$, $(\bE(\sigma_n)-\bE(\sigma^*))|_{\Omega_0}\in H^{1/2}(\Omega_0)^3$. Then it follows from Lemma \ref{regularity1} 
that 
$$\|\bE(\sigma_n)-\bE(\sigma^*)\|_{H^{1/2}(\Omega_0)^3}\leq C\|\sigma-\sigma^*\|_{L^2(\Omega_c)}\|\bJ_s\|_{L^2(\Omega)^3}.$$
From the convergence result in Lemma \ref{lem:3.1},  we know that $\bE(\sigma_n)|_{\Omega_0}\rightarrow\bE(\sigma^*)|_{\Omega_0}$ in $H^{1/2}(\Omega_0)^3$. Then by trace theorem, 
$\bn\times\bE(\sigma_n)|_\Gamma\rightarrow \bn\times\bE(\sigma^*)|_\Gamma \mbox{ in } L^2(\Gamma)^3$,
which concludes the compactness.
\finproof

\subsection{Non-uniqueness of the recovery of the conductivity}

In this subsection, we will study the uniqueness of the recovery of conductivity using the data 
$\bn\times\bE(\sigma)$ on the boundary 
$\Gamma$.  
Some techniques used here are motivated by the uniqueness argument in \cite{RCV}
for an inverse source problem.
Let $\bE_0$ be the background field which satisfies
\begin{eqnarray*}
\left\{
\begin{array}{rlll}
\nabla\times(\mu^{-1}\nabla\times\bE_0)-i\omega\sigma_0\bE_0&=&i\omega\bJ_s &\mbox{ in }\Omega,\\
\nabla\cdot\varepsilon\bE_0&=&0 &\mbox{ in }\Omega_0,
\end{array}
\right.
\end{eqnarray*}
with boundary conditions \eqref{BCs}.

If we have known the exact data $\bn\times\bE(\sigma)$ and the  background conductivity $\sigma_0$,  the recovery problem is reduced to determining $\sigma$ in $\Omega_c$ given $\bn\times(\bE(\sigma)-\bE_0)$ on $\Gamma$.
By simple calculation, we know that
\begin{eqnarray}
\left\{
\begin{array}{rlll}
\nabla\times(\mu^{-1}\nabla\times(\bE(\sigma)-\bE_0))-i\omega\sigma_0(\bE(\sigma)-\bE_0)&=&i\omega\sigma\bE(\sigma) &\mbox{ in }\Omega,\\
\nabla\cdot\varepsilon(\bE(\sigma)-\bE)&=&0 &\mbox{ in }\Omega_0,
\end{array}
\right.\label{pertb}
\end{eqnarray}
and $\bE(\sigma)-\bE_0$ satisfies boundary conditions \eqref{BCs}. 
Let $\bJ_e=i\omega\sigma\bE(\sigma)$. Then it is clear that $\bJ_e$ is supported in $\Omega_c$. In the following part of this section, we will consider an inverse source problem related to equation \eqref{pertb}. To be more specific, let $\bE$ satisfy the following equation
\begin{eqnarray}
\left\{
\begin{array}{rlll}
\nabla\times(\mu^{-1}\nabla\times\bE)-i\omega\sigma_0\bE&=&\bJ_e &\mbox{ in }\Omega,\\
\nabla\cdot\varepsilon\bE&=&0 &\mbox{ in }\Omega_0,
\end{array}
\right.\label{inversesource}
\end{eqnarray}
and boundary condition \eqref{BCs},  the corresponding inverse source problem is: 
\begin{eqnarray}
\mbox{ given data }\bn\times\bE\mbox{ on }\Gamma\mbox{, find the source }\bJ_e \mbox{ supported in }\Omega_c.\label{InvSource}
\end{eqnarray}
To proceed, we denote
\begin{eqnarray*}
W=\big\{\bu\in H(\curl
;\Omega_c)~\big|~\nabla\times(\mu^{-1}\nabla\times\bu)+i\omega\sigma_0\bu=0 \mbox{ in }\Omega_c,\bn\times\bu=0 \mbox{ on }\partial\Omega_c\setminus\Gamma_{0c}\big\}.
\end{eqnarray*}
It is easy to find that $W$ is not empty because the boundary value problem
\begin{eqnarray*}
\left\{
\begin{array}{clll}
\nabla\times(\mu^{-1}\nabla\times\bu)+i\omega\sigma_0\bu&=&0 &\mbox{ in }\Omega_c,\\
\bn\times\bu&=&0 &\mbox{ on }\partial\Omega_c\setminus\Gamma_{0c},\\
\bn\times\bu&=&\eta &\mbox{ on }\Gamma_{0c},
\end{array}
\right.
\end{eqnarray*}
has a unique solution for any $\eta\in H^{-1/2}(\Div;\Gamma_{0c})$.
Let $$L^2(\Omega_c)^3=W\oplus W^\bot.$$
We know that $W^\bot$ is not  trivial either. More precisely, for any $\bv\in C^\infty_0(\Omega_c)$, it is easy to find that $\nabla\times\mu^{-1}\nabla\times\bv-i\omega\sigma_0\bv\in W^\bot$.
The next theorem tells us there are non-radiation sources satisfying the inverse source problem \eqref{InvSource}.
\begin{thm}\label{NRsource}
 If $\bJ_e\in W^\bot$, the corresponding field is denoted by $\bE$, then $\bn\times\bE=0$ on $\Gamma_{0c}$ and $\Gamma$, in other words, $\bJ_e$ is a non-radiating source.
\end{thm}
\debproof
For any $\bu\in W$, by integration by parts, we have
\begin{eqnarray*}
&&0=\int_{\Omega_c}\bJ_e\cdot\overline\bu dx =\int_{\Omega_c}(\nabla\times(\mu^{-1}\nabla\times\bE)-i\omega\sigma_0\bE)\cdot\overline\bu dx\\
&&\ \ =\int_{\Omega_c}\bE\cdot(\nabla\times(\mu^{-1}\nabla\times\overline\bu)-i\omega\sigma_0\overline\bu)dx +\int_{\partial\Omega_c}\bn\times\mu^{-1}\nabla\times\bE\cdot\overline\bu+\mu^{-1}\bn\times\bE\cdot\nabla\times\overline\bu ds.
\end{eqnarray*}
Since $\bu\in W$, we obtain
\begin{eqnarray}
\int_{\Gamma_{0c}}\bn\times\mu^{-1}\nabla\times\bE\cdot\overline\bu+\mu^{-1}\bn\times\bE\cdot\nabla\times\overline\bu ds=0.
\label{nr-01}
\end{eqnarray}
For any $\eta\in H^{1/2}(\Div;\Gamma_{0c})$, let $\bw\in H_\Gamma(\curl;\Omega)$ be the solution of the following interface problem
\begin{eqnarray*}
\left\{
\begin{array}{rlll}
\nabla\times(\mu^{-1}\nabla\times\bw)+i\omega\sigma_0\bw&=&0 &\mbox{ in }\Omega_c\cup\Omega_0,\\
\nabla\cdot\varepsilon\bw&=&0 &\mbox{ in }\Omega_0,\\
\mu^{-1}\bn\times\nabla\times\bw|_{\Omega_c}&=&\mu^{-1}\bn\times\nabla\times\bw|_{\Omega_0}+\eta &\mbox{ on }\Gamma_{0c},
\end{array}
\right.
\end{eqnarray*}
and boundary conditions \eqref{BCs}.
With the similar method in Section \ref{sect2}, one can prove that the above system of equations is well-posed, i.e., for any $\eta\in H^{-1/2}(\Div;\Gamma_{0c})$, it has a unique solution $\bw$ in $H_\Gamma(\curl;\Omega)$. Furthermore, $\bw\neq 0$ in $\Omega_c$. If we choose $\bu=\bw|_{\Omega_c}$ in \eqref{nr-01}, it becomes
\begin{eqnarray*}
&&0= \int_{\Gamma_{0c}}\mu^{-1}\bn\times\nabla\times\bE\cdot\overline\bu+\mu^{-1}\bn\times\bE\cdot\nabla\times\overline\bu ds\\
&&\ \ =\int_{\Gamma_{0c}}\mu^{-1}\bn\times\nabla\times\bE\cdot\overline\bw ds-\int_{\Gamma_{0c}}\mu^{-1}\bn\times\nabla\times\overline\bw|_{\Omega_c}\cdot\bE ds\\
&&\ \ =\int_{\Gamma_{0c}}\mu^{-1}\bn\times\nabla\times\bE\cdot\overline\bw ds-\int_{\Gamma_{0c}}\mu^{-1}\bn\times\nabla\times\overline\bw|_{\Omega_0}\cdot\bE ds-\int_{\Gamma_{0c}}\overline\eta\cdot\bE ds.
\end{eqnarray*}
Using the fact that $\bn\times\bE=0$ on $\Gamma_D$ and $\bn\times\nabla\times\bw=0$ on $\Gamma$, we have
\begin{eqnarray}
\int_{\Gamma_{0c}}\mu^{-1}\bn\times\nabla\times\bE\cdot\overline\bw ds -\int_{\Gamma_{0c}}\overline\eta\cdot\bE ds-\int_{\partial\Omega_0}\mu^{-1}\bn\times\nabla\times\overline\bw\cdot\bE ds=0 .\label{nr-02}
\end{eqnarray}
By integration by parts,
\begin{eqnarray*}
&&\int_{\partial\Omega_0}\mu^{-1}\bn\times\nabla\times\overline\bw\cdot\bE ds=\int_{\Omega_0}\nabla\times(\mu^{-1}\nabla\times\overline\bw\cdot\bE)-\mu^{-1}\nabla\times\overline\bw\cdot\nabla\times\bE dx\\
&&\ \ =-\int_{\Omega_0}\mu^{-1}\nabla\times\overline\bw\cdot\nabla\times\bE dx\\
&&\ \ =-\int_{\Omega_0}\overline\bw\cdot\nabla\times(\mu^{-1}\nabla\times\bE)  dx +\int_{\partial\Omega_0}\mu^{-1}\bn\times\nabla\times\bE\cdot\overline\bw ds\\
&&\ \ =\int_{\partial\Omega_0}\mu^{-1}\bn\times\nabla\times\bE\cdot\overline\bw ds.
\end{eqnarray*}
Substituting the above results into \eqref{nr-02}, we have$$\int_{\Gamma_{0c}}\eta\cdot\bE ds=0.$$
Then $\bn\times\bE\times\bn=0$ on $\Gamma_{0c}$, and this implies that $\bE=0$ in $\Omega_0$ and $\bn\times\bE=0$ on $\Gamma$.
\finproof

With the help of Theorem \ref{NRsource}, we can give the following theorem about the non-uniqueness recovery property for the inverse eddy current problem.
 \begin{thm}\label{thm:uniq}
With the measurement satisfying $\bn\times(\bE(\sigma)-\bE_0)\neq 0$ on $\Gamma$, one can not determine $\sigma\bE(\sigma)$ uniquely.
 \end{thm}
 \debproof
Since $\bn\times(\bE(\sigma)-\bE)\neq 0$ on $\Gamma$, we can conclude that $\sigma\bE(\sigma)\notin W^\bot$ by Theorem \ref{NRsource}. Note that $\sigma\bE(\sigma)\in H_0(\curl;\Omega_c)$. If $\sigma\bE(\sigma)\in W$, we can conclude that $\sigma\bE(\sigma)$ is an homogeneous eigenfunction corresponding to imaginary eigenvalue $i\omega\mu\sigma_0$ in $\Omega_c$. It is impossible because the operator  $\nabla\times\nabla\times$ is a semi-positive operator on space $H_0(\curl;\Omega_c)$. So $\sigma\bE(\sigma)\notin W$ and we conclude that
$$i\omega\sigma\bE(\sigma)=\bJ_1+\bJ_2, \mbox{ and }\bJ_1\in W, \,\bJ_2\in W^\bot,\, \bJ_i\neq 0, \,i=1,2.$$
From $\bJ_2\neq 0$, by Theorem \ref{NRsource}, we finish the proof.
\finproof
\begin{rem}
We do not know whether the sources belong to $W$ can be uniquely determined or not. It does not matter because we know from Theorem \ref{thm:uniq} that there is always a non-radiation part of $\sigma\bE(\sigma)$.
\end{rem}

The secondary source term $\bJ_e=i\omega\bE(\sigma)$ depends on $\sigma$ nonlinearly. The following corollary tells us that when $\sigma$ is small enough, $\bJ_e$ can be approximated by $i\omega\bE_0$. 
Then $\bJ_e$ depends on $\sigma$ approximately and linearly.
\begin{cor}\label{cor3_1}
Assuming that $\sigma_0$ is a constant and $\sigma$ is small enough such that, for some $\sigma_m>0$,  $\sigma<\sigma_m<\sigma^2_0$, we have
\begin{eqnarray*}
\|\bE(\sigma)-\bE_0\|_{L^2(\Omega_c)^3}\leq \frac{\sigma_m}{\sigma^2_0-\sigma_m}\|\bE_0\|_{L^2(\Omega_c)^3}.
\end{eqnarray*}
\end{cor}
\debproof
Let $\hat\bE=\bE(\sigma)-\bE_0$. Multiplying both sides of the first equation in \eqref{pertb} by any $\bF\in H_\Gamma(\curl;\Omega)$, and using integration by parts, we have
\begin{eqnarray*}
\int_{\Omega}\mu^{-1}\nabla\times\hat\bE\cdot\nabla\times\overline\bF -i\omega\sigma_0\hat\bE\cdot\overline{\bF} dx=i\omega\int_{\Omega_c}\sigma\bE(\sigma)\cdot\overline\bF dx.
\end{eqnarray*}
Let $\bF=\hat\bE$, then
\begin{eqnarray*}
\|\mu^{-1/2}\nabla\times\hat\bE\|^2_{L^2(\Omega)}-i\omega\|\hat\bE\|^2_{L^2(\Omega_c)}=i\omega\int_{\Omega_c}\sigma\bE^{obs}\cdot\overline{\hat\bE} dx.
\end{eqnarray*}
So
\begin{eqnarray*}
\left\{
\begin{array}{rll}
\|\sigma_0\hat\bE\|^2_{L^2(\Omega_c)} &=& -\re\int_{\Omega_c}\sigma\bE(\sigma)\cdot\overline{\hat\bE}dx,\\
\|\mu^{-1/2}\nabla\times\hat\bE\|^2_{L^2(\Omega)}&=&\im \int_{\Omega_c}\sigma\bE(\sigma)\cdot\overline{\hat\bE} dx.
\end{array}
\right.
\end{eqnarray*}
From the first equality above, we can finish the proof with the following inequality
\begin{eqnarray*}
\|\sigma_0\hat\bE\|_{L^2(\Omega_c)^3}\leq \frac{\sigma}{\sigma_0}\|\bE^{obs}\|_{L^2(\Omega_c)^3}\leq \frac{\sigma_m}{\sigma_0}(\|\bE_0\|_{L^2(\Omega_c)^3}+\|\hat\bE\|_{L^2(\Omega_c)^3}).
\end{eqnarray*}
\finproof
\begin{rem}\label{rem:source}
We know that the secondary source $\bJ_e=i\omega\mu\sigma\bE(\sigma)=i\omega\mu\sigma\bE_0+i\omega\mu\sigma\hat\bE$. When $\sigma$ is small enough, 
Corollary \ref{cor3_1}
confirms that the second term is of high order of $\sigma$. If we drop the higher order term in the right-hand side of \eqref{pertb}, and assume that we can uniquely determine the secondary source, then with $\bE_0$ known we can uniquely determine $\sigma$, except that $\bE_0$ vanishes.
Unfortunately, with the same reason as in Theorem \ref{thm:uniq}, 
we know that 
$\sigma\bE_0$ does not lie in either $W$ or $W^\bot$, so we 
can not determine $\sigma\bE_0$ completely with the measurement on $\Gamma$.
\end{rem}

\subsection{Regularized inverse problem}
Since the eddy current inverse problem is ill-posed, we take the following regularization 
to transform the ill-posed problem to a problem that is at least mathematically stable with respect to the change 
of the noisy data for numerical solutions:
\begin{equation} 
\min_{\sigma\in H^1_0(\Omega_c)}\Phi_\alpha(\sigma):=\frac{1}{2}\|\bn\times(\bE(\sigma)-\bE^{obs})\|^2_{L^2(\Gamma)}+\frac{\alpha}{2}\|\nabla\sigma\|^2_{L^2(\Omega_c)},\label{opt_reg}
\end{equation}
where $\bE(\sigma)$ satisfies the constrained equation \eqref{MTE} or \eqref{MTE-w}, and $\alpha$ is the regularization parameter.

We first show the existence of the minimizers of the functional $\Phi_\alpha(\sigma)$. 

\begin{thm}\label{theorem:existence}
Under the same assumptions as in Theorem \ref{regularity}, there exists a minimizer $\sigma_\alpha$ to $\Phi_\alpha(\sigma)$ in $H^1_0(\Omega_c)$.
\end{thm}
\debproof
The proof is quite standard; 
see, e.g., \cite{CZ,FJZ} for the inverse elliptic and Maxwell problems.
But for readers's convenience, we give an outline of the proof, showing the main differences 
for the current eddy current problem.
First, we assume that $\{\sigma_n\}$ is a minimizing sequence for $\Phi_\alpha(\sigma)$, i.e.,
\begin{eqnarray*}
\lim_{n\rightarrow\infty}\Phi_\alpha(\sigma_n)=\inf_{\sigma\in H^1_0(\Omega_c)}\Phi_\alpha(\sigma).
\end{eqnarray*}
By  the convergence of $\{\Phi_\alpha(\sigma_n)\}$, we know that $\{\Phi_\alpha(\sigma_n)\}$ is bounded, 
so is $\{\|\nabla\sigma_n\|_{L^2(\Omega_c)}\}$. Therefore $\{\sigma_n\}$ is bounded in $H^1_0(\Omega_c)$, 
and there is a subsequence, still denoted by $\{\sigma_n\}$, converges weakly to $\sigma_\alpha$. 
This weak convergence is actually strong, due to the compact imbedding of $H^1(\Omega_c)$ 
in $L^2(\Omega_c)$. 
For $\sigma_n$ and $\sigma_\alpha$, we denote by $(\bE_n,\phi_n)$ and $(\bE_\alpha,\phi_\alpha)$
the solutions to the following two systems 
\begin{equation*}
\left\{
\begin{array}{clll}
a_n(\bE_n,\bF)+b(\nabla\phi_n,\bF) &=&i\omega\int_{\Omega} \bJ_s\cdot\overline\bF dx &~~\forall\,\bF\in H_\Gamma(\curl;\Omega), \\
b(\bE_n,\nabla\psi)&=&0&~~\forall\, \psi\in H^1_\Gamma(\Omega_0),
\end{array}
\right.
\end{equation*}
%
\begin{equation*}
\left\{
\begin{array}{clll}
a_\alpha(\bE_\alpha,\bF)+b(\nabla\phi_\alpha,\bF) &=&i\omega\int_{\Omega} \bJ_s\cdot\overline\bF dx&~~\forall\,\bF\in H_\Gamma(\curl;\Omega), \\
b(\bE_\alpha,\nabla\psi)&=&0&~~\forall\, \psi\in H^1_\Gamma(\Omega_0),
\end{array}
\right.
\end{equation*}
where $a_n(\bE,\bF)$ and $a_\alpha(\bE,\bF)$ are defined by 
\eqref{eq:a(u,v)} but with $\sigma$ there replaced by 
$\sigma_n$ and $\sigma_\alpha$, respectively. 
By Lemmas  \ref{regularity1} and \ref{lem:3.1}, the well-posedness of \eqref{MTE-w} and 
regularity results in Theorem \ref{regularity}, we have 
\begin{eqnarray*}
&&\bE_n|_{\Omega_0}\rightarrow \bE_\alpha|_{\Omega_0} \mbox{ ~in~ } H^{1/2}(\Omega_0)\,;
\quad 
\phi_n \rightarrow \phi_\alpha \mbox{ ~in~ } H^1_\Gamma(\Omega_0).
\end{eqnarray*}
Then by the trace theorem, $\bn\times\bE_n\rightarrow\bn\times\bE_\alpha$ in $L^2(\Gamma)$.
Using the strong convergence of $\bn\times\bE_n$, and the weakly lower semi-continuity 
of $\Phi_\alpha(\sigma)$, we get 
\begin{eqnarray*}
\Phi_\alpha(\sigma_\alpha)\leq \liminf_{n\rightarrow\infty}\Phi_\alpha(\sigma_n)=\inf_{\sigma\in H^1_0(\Omega_c)}\Phi_\alpha(\sigma)\,. 
\end{eqnarray*} 
\finproof

We end this section with the stability of the regularized optimization system \eqref{opt_reg}, 
whose proof can be done by standard arguments, 
along with some special techniques in the proof of Theorem\,\ref{theorem:existence}; 
see, e.g., \cite{CZ,FJZ} for the inverse elliptic and Maxwell problems.
\begin{thm}
Let $\{\bE_n\}$ be a sequence such that $\|\bn\times\bE_n-\bn\times\bE^{obs}\|_{L^2(\Gamma)}\rightarrow 0$ as $n\rightarrow \infty$, $\sigma_n$ be the minimizer of $\Phi_\alpha$ defined by \eqref{opt_reg} 
but with the quantity $\bn\times\bE^{obs}$ replaced by $\bn\times\bE_n$, 
then $\{\sigma_n\}$ has a subsequence which converges strongly 
to a minimizer of $\Phi_\alpha$ in $L^2(\Omega_c)$.
\end{thm}
\section{Nonlinear conjugate gradient method}\label{sect4}
In this section, we first introduce the Lagrangian of the optimization problem \eqref{opt_reg}, 
then derive the gradient of the objective functional and the G\^ateaux derivative of the electric field with respect to parameter $\sigma$. 
Finally, we formulate a nonlinear conjugate gradient method with an approximate scheme for 
step length.

\subsection{Lagrangian for the continuous optimization problem}
In order to calculate the gradient of the objective functional $\Phi_\alpha$ with respect to $\sigma$, we use the standard adjoint state technique. We first recast the problem \eqref{opt_reg} into an unconstrained optimization 
by introducing some multipliers to relax the PDE constraint. Since the system of equations \eqref{MTE} and its weak formulation \eqref{MTE-w} are complex-valued, we relax the constraint in the real and imaginary parts 
separately to reformulate them into a real-valued unconstrained optimization. 
Let $\bE=\rE_1+i\rE_2$, $\phi=\phi_1+i\phi_2$, $ i\omega\mu\bJ_s=\bbf_1+i\bbf_2$ and define $a_i:H_\Gamma(\curl;\Omega)\times H_\Gamma(\curl;\Omega)\rightarrow\mathbb{R}$ for $i=1,2$ as
\begin{eqnarray}
a_1(\bE, {\rm F})=\int_{\Omega}\mu^{-1}\nabla\times\rE_1\cdot\nabla\times {\rm F} dx +\int_{\Omega_c}\omega(\sigma_0+\sigma)\rE_2\cdot {\rm F} dx, \label{eq:a1}\\
a_2(\bE, {\rm F})=\int_{\Omega}\mu^{-1}\nabla\times\rE_2\cdot\nabla\times{\rm F} dx -\int_{\Omega_c}\omega(\sigma_0+\sigma)\rE_1\cdot {\rm F} dx. \label{eq:a2}
\end{eqnarray}
By taking the test functions in real function spaces $H_\Gamma(\curl;\Omega)$ and $H^1_\Gamma(\Omega_0)$,
the complex-valued system \eqref{MTE-w} becomes the following real-valued system for $i=1,2$:
\begin{eqnarray}
\left\{\begin{array}{clll}
a_i(\bE,\rF_i)+b(\nabla\phi_i,\rF_i) &=&\int_{\Omega}\bbf_i\cdot\rF_i dx &~~\forall\, \rF_i\in H_\Gamma(\curl; \Omega),\\
b(\rE_i,\psi_i)&=&0 & ~~\forall\,\psi_i\in H^1_\Gamma(\Omega_0)\,.
\end{array}\right.\label{state0}
\end{eqnarray}

Accordingly, we rewrite $\Phi_\alpha(\sigma)$ in \eqref{opt_reg} as $\Phi_\alpha(\bE,\sigma)$, that is,
 \begin{eqnarray*}
 \Phi_\alpha(\bE,\sigma)=
 \frac{1}{2}\big(\|\bn\times(\rE_1-\rE^{obs}_1)\|^2_{L^2(\Gamma)}+\|\bn\times(\rE_2-\rE^{obs}_2)\|^2_{L^2(\Gamma)}\big)
 +\frac{\alpha}{2}\|\nabla\sigma\|^2_{L^2(\Omega_c)},
 \end{eqnarray*}
 where $\rE^{obs}_1$ and $\rE^{obs}_2$ are the real and imaginary parts of $\bE$, respectively.
 Now we use $\Sigma$ to denote the product space $H_\Gamma(\curl,\Omega)\times H_\Gamma(\curl,\Omega)\times H^1_\Gamma(\Omega_0)\times H^1_\Gamma(\Omega_0)$, then we can define a Lagrangian functional $L$ from $\Sigma\times \Sigma\times H^1_0(\Omega_c)$ to $\mathbb{R}$:
\begin{eqnarray}
&&L((\rE_1,\rE_2,\phi_1,\phi_2),(\rF_1,\rF_2,\psi_1,\psi_2),\sigma)\nonumber\\
&=&\Phi_\alpha(\bE,\sigma)+\sum_{i=1}^2\big(a_i(\bE,\rF_i)+b(\nabla\phi_i,\rF_i)-\int_\Omega\bbf_1\cdot\rF_idx+b(\rE_i,\psi_i)\big)\label{lagrangian},
\end{eqnarray}
where real functions $\rF_1,\rF_2,\psi_1,\psi_2$ are Lagrange multipliers. Using the adjoint state technique, 
we can deduce that
\begin{eqnarray}
\frac{\partial L}{\partial \sigma}(\tilde\sigma)=\alpha\int_{\Omega_c}\nabla\sigma\cdot\nabla\tilde\sigma dx -\omega\int_{\Omega_c}(\rE_1\cdot\rF_2-\rE_2\cdot\rF_1)\tilde\sigma dx ~~\forall\,\tilde\sigma\in H^1_0(\Omega_c), \label{grad_w0}
\end{eqnarray}
where $\rE_1,\rE_2$ and $\rF_1, \rF_2$ are the solutions of the following systems
\begin{eqnarray}
\frac{\partial L}{\partial(\rF_1,\rF_2,\psi_1,\psi_2)}((\tilde{\rF}_1,\tilde{\rF}_2,\tilde{\psi}_1,\tilde{\psi}_2))=0 & ~~\forall\, (\tilde{\rF}_1,\tilde{\rF}_2,\tilde{\psi}_1,\tilde{\psi}_2)\in \Sigma \label{Lg1},\\
\frac{\partial L}{\partial(\rE_1,\rE_2,\phi_1,\phi_2)}((\tilde{\rE}_1,\tilde{\rE}_2,\tilde{\phi}_1,\tilde{\phi}_2))=0 & ~~\forall\, (\tilde{\rE}_1,\tilde{\rE}_2,\tilde{\phi}_1,\tilde{\phi}_2)\in \Sigma \label{Lg2}.
\end{eqnarray}
We can easily check \eqref{Lg1} is exactly the state system \eqref{state0}, 
while \eqref{Lg2}  yields its adjoint-state system
\begin{eqnarray}
\left\{\begin{array}{rlll}
a_j(\hat\bF,\tilde\rE_i)+b(\nabla\psi_i,\tilde\rE_i) &=&\int_\Gamma\Delta\rE_i\cdot\tilde\rE_i ds &~~\forall\, \tilde\rE_i\in H_\Gamma(\curl,\Omega), ~~j=2,1 \mbox{ for } i=1,2,\\
b(\rF_i,\tilde\phi_i)&=&0 & ~~\forall\,\tilde\phi_i\in H^1_\Gamma(\Omega_0),  ~~i=1,2,
\end{array}\right. \label{adjoint0}
\end{eqnarray}
where
$\Delta\rE_i=\bn\times(\rE_i^{obs}-\rE_i)\times\bn$ for $i=1,2$,  and $\hat\bF=\rF_2+i\rF_1$.

Let $(\bE,\phi_1,\phi_2)$ be the solution of system \eqref{state0}  and $\phi=\phi_1+i\phi_2$, 
then $(\bE,\phi)$ solves the system
\begin{eqnarray}
\left\{
\begin{array}{rlll}
a(\bE,\tilde\bF)+ b(\nabla\phi,\tilde\bF)&=&\int_{\Omega} (\bbf_1+i\bbf_2)\cdot\overline{\tilde\bF} dx&~~\forall\, \tilde{\bF}\in H_\Gamma(\curl,\Omega),\\
b(\bE,\nabla\tilde\psi)&=&0 &~~\forall\,\tilde{\psi}\in  H^1_\Gamma(\Omega_0).
\end{array}
\right. \label{state-w}
\end{eqnarray}
It is easy to find that equation \eqref{state-w} is equivalent to \eqref{MTE-w} because $\bbf_1+i\bbf_2=i\omega\mu\bJ_s$.

On the other hand, let $(\hat\bF,\psi_1,\psi_2)$ be the solution of system \eqref{adjoint0},
and let $\bF=-i\hat\bF=\bF_1-i\bF_2$, $\psi=\psi_1-i\psi_2$. Then we can check that $(\bF,\psi)$ solves 
the system 
\begin{eqnarray}
\left\{
\begin{array}{rlll}
a(\bF,\tilde\bE)+ b(\nabla\psi,\tilde\bE)&=&\int_\Gamma\bn\times\overline{(\bE^{obs}-\bE)}\times\bn\cdot\overline{\tilde\bE} ds &~~\forall\, \tilde\bE\in H_\Gamma(\curl,\Omega), \\
b(\bF,\nabla\tilde\phi)&=&0 &~~\forall\,\tilde{\phi}\in  H^1_\Gamma(\Omega_0).
\end{array}
\right. \label{adjoint-w}
\end{eqnarray}

Now we show that the relation \eqref{grad_w0} gives the gradient of $\Phi_\alpha(\sigma)$ 
in the weak sense. To see this, we define $g(\sigma)\in L^2(\Omega_c)$ with
\begin{eqnarray}
\int_{\Omega_c}g(\sigma)\tilde\sigma dx = 
\alpha\int_{\Omega_c}\nabla\sigma\cdot\nabla\tilde\sigma dx +\omega\int_{\Omega_c}\im{(\bE\cdot\bF)}\tilde\sigma dx ~~\forall\,\tilde\sigma\in H^1_0(\Omega_c), \label{grad_w}
\end{eqnarray}
then we have in the weak sense that 
\begin{eqnarray}
\frac{\partial \Phi_\alpha(\sigma)}{\partial\sigma} = g(\sigma) .\label{grad}
\end{eqnarray}
We can solve the equations \eqref{state-w}, \eqref{adjoint-w} and \eqref{grad_w} to calculate the gradient of the objective functional $\Phi_\alpha$ with respect to $\sigma$.

\subsection{Adjoint-state equations}
The adjoint system of equations \eqref{adjoint-w} looks similar to the state system \eqref{state-w} formally, but they are quite different in terms of their corresponding differential equations, which 
we derive below explicitly and explain their main differences. 
To do so, for the solution $(\bF,\psi)$ to \eqref{adjoint-w}, 
we can derive by integration by parts that
\begin{eqnarray*}
&&\int_{\Omega}(\nabla\times(\mu^{-1}\nabla\times\bF)-i\omega(\sigma_0+\sigma)\bF)\cdot\tilde\bE dx+\int_{\Omega_0}\varepsilon\nabla\psi\cdot\tilde\bE dx
-\int_\Gamma\mu^{-1}\bn\times\nabla\times\bF\cdot\tilde\bE ds\\
&&\ \ \ =\int_\Gamma(\bn\times\overline{(\bE^{obs}-\bE)}\times\bn)\cdot\tilde\bE ds.
\end{eqnarray*}
On the other hand, for any function $\phi\in H^1_\Gamma(\Omega_0)$, we extend it to $\Omega$ by zero, 
then choosing $\tilde\bE=\nabla\phi$ in \eqref{adjoint-w}, we obtain 
\begin{eqnarray*}
\int_{\Omega_0}\varepsilon\nabla\psi\cdot\nabla\phi dx=\int_\Gamma\bn\times\overline{(\bE^{obs}-\bE)}\times\bn\cdot\nabla\phi ds.
\end{eqnarray*}
This gives the corresponding differential equation for $\psi$:
\begin{eqnarray}
\left\{
\begin{array}{clll}
\nabla\cdot(\varepsilon\nabla\psi)&=&0 &\mbox{ in }\Omega_0,\\
\varepsilon\frac{\partial\psi}{\partial\bn}&=&\Div_\tau(\bn\times\overline{(\bE^{obs}-\bE)}\times\bn)&\mbox{ on }\Gamma,\\
\psi&=&0 &\mbox{ on }\partial\Omega_0\setminus\Gamma.
\end{array}
\right.\label{psid}
\end{eqnarray}
If we choose $\tilde\bE\in H_\Gamma(\curl,\Omega)$ and $\bn\times\tilde{\bE}=0$ on $\Gamma$, then $\bF$ needs to satisfy
\begin{eqnarray*}
\nabla\times(\mu^{-1}\nabla\times\bF)-i\omega(\sigma_0+\sigma)\bF+\varepsilon\psi=0 \mbox{ in } \Omega.
\end{eqnarray*}
If we choose $\bn\times\tilde\bE\neq 0$ on $\Gamma$, we can derive the following boundary condition that $\bF$ need to satisfy:
\begin{eqnarray*}
\mu^{-1}\bn\times\nabla\times\bF=-\bn\times\overline{(\bE^{obs}-\bE)}\times\bn \mbox{ on }\Gamma.
\end{eqnarray*}
Together with the  second equation in \eqref{adjoint-w}, we have derived 
the system of differential equations for the solution $\bF$ to  \eqref{adjoint-w}:
\begin{eqnarray}
\left\{
\begin{array}{rlll}
\nabla\times(\mu^{-1}\nabla\times\bF)-i\omega(\sigma_0+\sigma)\bF&=&-\varepsilon\nabla\psi &\mbox{ in } \Omega,\\
\nabla\cdot\varepsilon\bF &=&0 &\mbox{ in }\Omega_0, \\
\mu^{-1}\bn\times\nabla\times\bF&=&-\bn\times\overline{(\bE^{obs}-\bE)}\times\bn & \mbox{ on }\Gamma,\\
\bn\cdot\bF&=&0 &\mbox{ on }\Gamma,\\
\bn\times\bF &=&0 &\mbox{ on }\Gamma_D.
\end{array}
\right.\label{ajdiff}
\end{eqnarray}
This, along with \eqref{psid}, provides the differential equations of the solution $(\bF,\psi)$ 
to equation \eqref{adjoint-w}. 

Now we study a special case when $\Div_\tau(\bn\times(\bE^{obs}-\bE)\times\bn)=0$ on $\Gamma$. 
Then we know $\psi=0$ from \eqref{psid}, and equation \eqref{ajdiff} reduces to 
\begin{eqnarray}
\left\{
\begin{array}{rlll}
\nabla\times(\mu^{-1}\nabla\times\bF)-i\omega(\sigma_0+\sigma)\bF&=&0 &\mbox{ in } \Omega,\\
\nabla\cdot\varepsilon\bF &=&0 &\mbox{ in }\Omega_0, \\
\mu^{-1}\bn\times\nabla\times\bF&=&-\bn\times\overline{(\bE^{obs}-\bE)}\times\bn & \mbox{ on }\Gamma,\\
\bn\cdot\bF&=&0 &\mbox{ on }\Gamma,\\
\bn\times\bF &=&0 &\mbox{ on }\Gamma_D.
\end{array}
\right.\label{ajdiff1}
\end{eqnarray}
By the definition of the surface divergence \cite{BCS, Monk}, we have
\begin{eqnarray}\label{eq:div_tau}
\Div_\tau(\bn\times(\bE^{obs}-\bE)\times\bn)=-\bn\cdot\nabla\times((\bE^{obs}-\bE)\times\bn)|_\Gamma.
\end{eqnarray}
Using the relation that $\nabla\times(\bu\times\bv)=\bu(\nabla\cdot\bv)-(\bu\cdot\nabla)\bv+(\bv\cdot\nabla)\bu-\bv(\nabla\cdot\bu) \forall\, \bu,\bv$ and the normal $\bn=(0,0,1)$ of $\Gamma$ in our case, we get 
\begin{eqnarray*}
\nabla\times((\bE^{obs}-\bE)\times\bn)=(\bn\cdot\nabla)(\bE^{obs}-\bE)-\bn(\nabla\cdot(\bE^{obs}-\bE)).
\end{eqnarray*}
Assuming $\bE^{obs}$ is the solution to system \eqref{MTE} with the true conductivity, 
we then have $\nabla\cdot(\bE^{obs}-\bE)=0$ in $\Omega_0$. This leads to
\begin{eqnarray*}
\nabla\times((\bE^{obs}-\bE)\times\bn)=(\bn\cdot\nabla)(\bE^{obs}-\bE).
\end{eqnarray*}
Therefore we deduce from \eqref{eq:div_tau} that 
\begin{eqnarray*}
\Div_\tau(\bn\times(\bE^{obs}-\bE)\times\bn)&=&-\bn\cdot(\bn\cdot\nabla)(\bE^{obs}-\bE)|_\Gamma\\
&=&-(\bn\cdot\nabla) (\bE^{obs}-\bE)\cdot\bn|_\Gamma\\
&=&-\frac{\partial((\bE^{obs}-\bE)\cdot\bn)}{\partial\bn}.
\end{eqnarray*}
We know that the above relation is valid for any constant vector $\bn$. For the current case 
with $\bn=(0,0,1)$ and $\varepsilon$ being a constant on $\Gamma$, the above derivation can be simplified. Let $\bE_x,\bE_y,\bE_z$ be the components of $\bE$ alone $x,y,z$-axis, respectively, then 
\begin{eqnarray*}
\Div_\tau (\bn\times(\bE^{obs}-\bE)\times\bn)&=&\frac{\partial(\bE^{obs}-\bE)_x}{\partial x}+\frac{\partial(\bE^{obs}-\bE)_y}{\partial y}
=-\frac{\partial(\bE^{obs}-\bE)_z}{\partial z}.
\end{eqnarray*}
So the condition that $\Div_\tau (\bn\times(\bE^{obs}-\bE)\times\bn)=0$ is equivalent to
\begin{eqnarray*}
\frac{\partial(\bE^{obs}-\bE)_z}{\partial z}=0 \mbox{ on }\Gamma.
\end{eqnarray*}

In general, the condition that $\Div_\tau (\bn\times(\bE^{obs}-\bE)\times\bn)=0$ 
is not true, so we do not have $\psi=0$.
Comparing with \eqref{MTE}, we can see that the adjoint equation has a special source 
$\varepsilon\nabla\psi$, where $\psi$ solves the equation \eqref{psid}. 
%
Provided that $\Div_\tau(\bn\times\overline{(\bE^{obs}-\bE)}\times\bn)\in H^{-1/2}(\Gamma)$,
the equation \eqref{psid} is well-posed, hence 
the adjoint system \eqref{adjoint-w} is well-posed, due to the well-posedness of \eqref{MTE-w}. 

\begin{rem}
Generally speaking, if $\bE^{obs}-\bE \in H_\Gamma(\curl;\Omega)$, we have that $\bn\times(\bE^{obs}-\bE)\times\bn \in H^{-1/2}(\rm{Curl};\Gamma)$, the dual space of $H^{-1/2}(\Div;\Gamma)$ \cite{BCS}. Then $\Div_\tau(\bn\times(\bE^{obs}-\bE)\times\bn)$ may not belong to $H^{-1/2}(\Gamma)$. But with the discussion of 
the regularity in Section \ref{sect2} and the fact that $\Gamma$ is part of the boundary of a convex domain, 
we have that $\bn\times(\bE^{obs}-\bE)\times\bn \in L^2(\Gamma)$ 
and $\Div_\tau(\bn\times(\bE^{obs}-\bE)\times\bn) \in H^{-1/2}(\Gamma)$.
\end{rem}

\subsection{G\^ateaux derivative of the electric field $\bE$}
In this subsection we derive the G\^ateaux derivative of the electric field with respect to the conductivity
$\sigma$. The derivative is needed to compute at each iteration of 
the nonlinear conjugate gradient algorithm (cf.\,Subsection \ref{nlcgalg}). 

For any 
 $\sigma\in H^1_0(\Omega_c)$,  we write $\sigma=\sigma_a+\sigma_b$ with $\sigma_a,\sigma_b\in H^1_0(\Omega_c)$, and decompose the corresponding 
solution $\bE(\sigma)$ to the system \eqref{MTE} as $\bE=\bE_0+\bE_1+\bE_2$, 
where $\bE_0:=\bE_0(\sigma_0+\sigma_a)$, $\bE_1:=\bE_1(\sigma_0+\sigma_a;\sigma_b)$ 
and $\bE_2:=\bE_2(\sigma_0+\sigma; \sigma_b)$ solve the following systems, respectively,  
along with boundary conditions \eqref{BCs}, 
\begin{eqnarray*}
\left\{\begin{array}{rlll}
\nabla\times(\mu^{-1}\nabla\times \bE_0)-i\omega(\sigma_0+\sigma_a) \bE_0&=&i\omega \bJ_s & \mbox{ in } \Omega, \\
\nabla\cdot\varepsilon \bE_0&=&0 &\mbox{ in }\Omega_0, 
\end{array}
\right.
\end{eqnarray*} 
\begin{eqnarray}
\left\{\begin{array}{rlll}
\nabla\times(\mu^{-1}\nabla\times \bE_1)-i\omega(\sigma_0+\sigma_a) \bE_1&=&i\omega\mu\sigma_b\bE_0 & \mbox{ in } \Omega, \\
\nabla\cdot\varepsilon \bE_1&=&0 &\mbox{ in }\Omega_0,  
\end{array}
\right.\label{gateaux}
\end{eqnarray} 
\begin{eqnarray*}
\left\{\begin{array}{rlll}
\nabla\times(\mu^{-1}\nabla\times \bE_2)-i\omega(\sigma_0+\sigma) \bE_2
&=&i\omega\mu\sigma_b\bE_1 & \mbox{ in } \Omega, \\
\nabla\cdot\varepsilon \bE_2&=&0 &\mbox{ in }\Omega_0\,.   
\end{array}
\right.
\end{eqnarray*} 
With the help of Corollary \ref{cor3_1}, we know that $\|\bE_1(\sigma_0+\sigma_a;\sigma_b)\|_{H(\curl;\Omega)}\leq C\|\sigma_b\|_{L^2(\Omega_c)}$ for small $\sigma_b$, 
and a simple integration by parts gives the following estimate of $\bE_2$:
\begin{eqnarray*}
\|\bE_2\|_{H(\curl;\Omega)}\leq C\|\sigma_b\|_{L^2(\Omega_c)}\|\bE_1\|_{L^2(\Omega)}\leq C\|\sigma_b\|^2_{L^2(\Omega_c)}.
\end{eqnarray*} 
This leads to
\begin{eqnarray*}
\lim_{\|\sigma\|_{L^2(\Omega_c)}\rightarrow 0}\frac{\|\bE(\sigma_0+\sigma)-\bE_0(\sigma_0+\sigma_a)-\bE_1(\sigma_0+\sigma_a;\sigma_b)\|_{H(\curl;\Omega)}}{\|\sigma\|_{L^2(\Omega_c)}}=0\,,
\end{eqnarray*} 
hence we know $\bE_1(\sigma_0+\sigma_a;\sigma_b)$ gives the G\^ateaux derivative of $\bE$ along the direction $\sigma_b$ at $\sigma_0+\sigma_a$. Since $\bE_1$ depends on $\sigma_b$ linearly, we have for given $\sigma_b$ and small $\gamma$ that 
\begin{eqnarray}
\bE(\sigma_0+\sigma_a+\gamma\sigma_b)=\bE_0(\sigma_0+\sigma_a)+\gamma\bE_1(\sigma_0+\sigma_a;\sigma_b)+o(|\gamma|\|\sigma_b\|_{L^2(\Omega_c)}). \label{linapp}
\end{eqnarray}
We note that the first two terms in the right-hand side above 
is the linear approximation of the electric field $\bE(\sigma_0+\sigma_a+\gamma\sigma_b)$.  With this approximation, let 
\begin{eqnarray}
 \Psi(\gamma)=
 \frac{1}{2}\|\bn\times(\bE_0+\gamma\bE_1-\bE^{obs})\|^2_{L^2(\Gamma)}
 +\frac{\alpha}{2}\|\nabla(\sigma_a+\gamma\sigma_b)\|^2_{L^2(\Omega_c)}, \label{linearap}
 \end{eqnarray}
 we have
 \begin{eqnarray*}
 \Phi_\alpha(\bE,\sigma_0+\sigma_a+\gamma\sigma_b)=\Psi(\gamma)+o(|\gamma|\|\sigma_b\|_{L^2(\Omega_c)}).
 \end{eqnarray*}
It is easy to find that $\Psi(\gamma)$ is a quadratic function with respect to $\gamma$, 
which we use to help us compute the descent step size in our iterative Algorithm \ref{NLCG}.

\subsection{Finite element discretization of the minimization problem}
In this section we discuss the edge element approximation of 
the optimization system \eqref{opt_reg}. For this purpose, we
partition the domain $\Omega$ into a set of tetrahedral elements $\mathcal{M}_h$, 
with each element $K\in\mathcal{M}_h$ lying completely in $\Omega_c$ or $\Omega_0$ 
Let $\mathcal{M}^0_h$ and $\mathcal{M}_h^c$ be the unions of elements contained in $\Omega_0$ and $\Omega_c$, respectively. 
Then we define the N\'ed\'elec edge element space 
  \begin{eqnarray*}
  \bX_h=\big\{\bu_h\in H_\Gamma(\curl;\Omega)~\big|~ \bu_h|_K=a_K+b_K\times\bx, a_K,b_K\in\mathbb{R}^3\big\}.
  \end{eqnarray*}
 and $U_h\subset H^1_\Gamma(\Omega_0)$ and $V_h\subset H^1_0(\Omega_c)$ 
 to be the standard continuous piecewise linear finite element spaces 
 over $\mathcal{M}^0_h$ and $\mathcal{M}^c_h$ respectively. For ease of presentation, we use the notation $\Sigma_h=\bX_h\times\bX_h\times U_h\times U_h$ in the sequel. 
With these preparations, we propose the approximation of the optimization \eqref{opt_reg}:
\begin{eqnarray}
\min_{\sigma_h\in V_h}\Phi_\alpha(\sigma_h)=\frac{1}{2}\|\bn\times(\bE_h(\sigma_h)-\bE^{obs})\|^2_{L^2(\Gamma)}
+\frac{\alpha}{2}\|\nabla\sigma_h\|^2_{L^2(\Omega_c)} \label{obj-discrete}
\end{eqnarray}
where $\bE_h(\sigma_h)$ solves 
\begin{eqnarray}
\left\{
\begin{array}{rlll}
a_h(\bE_h,\bF_h)+b(\nabla\phi_h,\bF_h)&=&i\omega\int_{\Omega} \bJ_s\cdot\overline{\bF_h} dx&~~\forall\, \bF_h\in \bX_h,\\
b(\bE_h,\nabla\psi_h)&=&0&~~\forall\, \psi_h\in U_h,
\end{array}
\right. \label{MTE-discrete}
\end{eqnarray}
where $a_h$ is given by the sesquilinear operator $a_h(\bE,\bF)=\int_{\Omega}\mu^{-1}\nabla\times\bE_h\cdot\nabla\times\overline{\bF_h}-i\omega(\sigma_0+\sigma_h)\bE_h\cdot\overline{\bF_h}dx$ for all $\bE$, $\bF\in \bX_h$.
By writing the space
$$\mathbf{Y}_h=\big\{\bu_h\in\bX_h~\big|~b(\bu_h,\nabla\phi_h)=0 ~~\forall\,\phi_h\in U_h\big\},$$
we know \eqref{MTE-discrete} is a saddle-point problem that is equivalent to the problem
$\bE_h\in\mathbf{Y}_h$ satisfying 
\begin{eqnarray}
a_h(\bE_h,\bF_h)=i\omega\int_\Omega\bJ_s\cdot\overline{\bF_h} dx~~\forall\,\bF_h\in\mathbf{Y}_h.\label{MTE0h}
\end{eqnarray}
We can easily see that $\mathbf{Y}_h$ is not a subspace of $\mathbf{Y}$, so we can not deduce 
the well-posedness of \eqref{MTE0h} from that of the continuous weak problem \eqref{MTE0}. 
Instead the well-posednss of \eqref{MTE0h} can be achieved from that 
of \eqref{MTE-w} by using the fact that $\Sigma_h$ is a subspace of $\Sigma$ and 
following the arguments in \cite{CCCZ} \cite{CXZ} for the magneto-static problem and 
field/circuit coupling problem. 
%

Similarly to the proof of Theorem\,\ref{theorem:existence}, we have the following existence. 
\begin{thm}
There exists at least one minimizer to the discrete optimization problem \eqref{obj-discrete}.
\end{thm}

Now we introduce a discrete Lagrangian on $\Sigma_h\times\Sigma_h\times V_h$ associated 
with \eqref{obj-discrete}. To do so, we first
define $a_{i,h}(\cdot,\cdot)$ for $i=1,2$ to be 
the same bilinear form as $a_i(\cdot,\cdot)$ defined in \eqref{eq:a1}-\eqref{eq:a2}, 
but with $\sigma$ replaced by $\sigma_h$. 
Then we define the discrete Lagrangian as 
\begin{eqnarray*}
&&L((\bE_h,\phi_h),(\bF_h,\psi_h),\sigma_h)\\
&&\ \ \ =\Phi_\alpha(\bE_h,\sigma_h)+\sum^2_{i=1}\big(a_{i,h}(\bE^i_h,\rF^i_h)+b(\nabla\phi^i_h,\rF^i_h)-\int_\Omega\bbf_1\cdot\rF^i_hdx+b(\rE^i_h,\psi^i_h)\big).
\end{eqnarray*}
Let $g_h(\sigma_h):=\frac{\partial \Phi_\alpha(\sigma_h)}{\sigma_h}$, 
then we can derive a similar relation to the continuous one \eqref{grad_w}:
\begin{eqnarray}
\int_{\Omega_c}g_h \tilde\sigma_h dx=
\alpha\int_{\Omega_c}\nabla\sigma_h\cdot\nabla\tilde\sigma_h dx -\omega\int_{\Omega_c}(\rE_h^1\cdot\rF_h^2-\rE_h^2\cdot\rF_h^1)\tilde\sigma_h dx ~~\forall\,\tilde\sigma_h\in V_h,\label{grad0d}
\end{eqnarray}
where $\bE_h=\rE_h^1+i\rE_h^2$ and $\hat\bF_h=\rF_h^2+i\rF_h^1$ solve the following state  
and adjoint systems, respectively, 
\begin{eqnarray}
\left\{\begin{array}{rlll}
a_{i,h}(\bE_h,\tilde\rF_h^i)+b(\nabla\phi_h^i,\tilde\rF_h^i) &=&\int_{\Omega}\bbf_1\cdot\tilde\rF_h^i dx &~~\forall\, \tilde\rF_h^i\in \bX_h,i=1,2,\\
b(\rE_h^i,\tilde\psi_h^i)&=&0 & ~~\forall\,\tilde\psi_h^i\in U_h, i=1,2,
\end{array}\right.\label{state0d}
\end{eqnarray}
\begin{eqnarray}
\left\{\begin{array}{rlll}
a_{j,h}(\hat\bF,\tilde\rE_h^i)+b(\nabla\psi_h^i,\tilde\rE_h^i) &=&\int_\Gamma\bn\times(\bE_i^{obs}-\rE_h^i)\times\bn\cdot\tilde\bE_h^1 ds &\forall\,\tilde\rE_h^i\in \bX_h, j=2,1 \mbox{ and } i=1,2,\\
b(\rF_h^i,\tilde\phi_h^i)&=&0 &\forall\,\tilde\phi_h^i\in U_h, i=1,2.
\end{array}\right. \label{adjoint0d}
\end{eqnarray}
If we write 
$\phi_h=\phi^1_h+i\phi^2_h$, $\bF_h=\bF^1_h-i\bF^2_h$, $\psi_h=\psi^1_h-i\psi^2_h$, 
then equations \eqref{state0d}-\eqref{adjoint0d} are just the discrete versions 
of \eqref{state0}-\eqref{adjoint0}, with their corresponding complex-valued systems given by 
\begin{eqnarray}
\left\{
\begin{array}{rlll}
a_h(\bE_h,\tilde\bF_h)+ b(\nabla\phi_h,\tilde\bF_h)&=&\int_{\Omega} (\bbf_1+i\bbf_2)\cdot\overline{\tilde\bF_h} dx&~~\forall\, \tilde{\bF}_h\in X_h,\\
b(\bE_h,\nabla\tilde\psi_h)&=&0&~~\forall\,\tilde{\psi}_h\in  U_h,
\end{array}
\right. \label{state-wd}
\end{eqnarray}
\begin{eqnarray}
\left\{
\begin{array}{rlll}
a_h(\bF_h,\tilde{\bE}_h)+ b(\nabla\psi_h,\tilde \bE_h)&=&\int_{\Gamma}\bn\times\overline{( \bE^{obs}-\bE_h)}\times\bn\cdot\overline{\tilde{\bE}_h} ds&~~\forall\, \tilde{\bE}_h\in X_h, \\
b(\bF_h,\nabla\tilde\phi_h)&=&0&~~\forall\,\tilde{\phi}_h\in  U_h.
\end{array}
\right. \label{adjoint-wd}
\end{eqnarray}
In addition, we can see that equation \eqref{grad0d} can be simplified as
\begin{eqnarray}
\int_{\Omega_c}g_h \tilde\sigma_h dx=\alpha\int_{\Omega_c}\nabla\sigma_h\cdot \nabla\tilde\sigma_h dx+\int_{\Omega_c} \omega\im(\bE_h\cdot\bF_h)\tilde\sigma_h dx ~~\forall\, \tilde{\sigma}_h\in V_h.\label{grad_wd}
\end{eqnarray}

\subsection{A nonlinear conjugate gradient method}\label{nlcgalg}
With the derivations in the previous subsections, we can now formulate the following 
nonlinear conjugate gradient algorithm for solving the discrete optimization problem \eqref{obj-discrete}.

\begin{alg}[NLCG method]\label{NLCG} 
Given the observation data $\bn\times\bE^{obs}$ on $\Gamma$,
the background medium $\sigma_0$, the initial guess $\sigma_h^0$;  set $k=0$.
\begin{itemize}
\item[1.] Solve problem \eqref{state-wd} with $\sigma_h= \sigma^k_h$ to get $\bE^k_h$ and $\phi^k_h$;
\item[2.] Solve problem \eqref{adjoint-wd} with $\sigma_h=\sigma^k_h$ to get $\bF^k_h$ and $\psi^k_h$;
\item[3.] Solve problem \eqref{grad_wd} to get the gradient $g^k_h$;
\item[4.] Update the descent direction
 $d_k=-g_h^k+\beta_k d_{k-1},$
with the step size $\beta_k$ computed by 
$$\beta_k=\left\{\begin{array}{cc}\frac{\|g_h^k\|^2_{L^2(\Omega_c)}}{\|g_h^{k-1}\|^2_{L^2(\Omega_c)}}\, 
&\mbox{for}~k>0,\\
 0\, &\mbox{for} ~k=0;\end{array}\right.$$
\item[5.] Solve problem \eqref{gateaux} with $\sigma_a=\sigma_h^k$ and $ \sigma_b= d_k$
for the solution $\bE_h^g$;
\item[6.] Compute 
$$\gamma_k=-\frac{\int_\Gamma\Re((\bE_h^k-\bE^{obs})\times\bn\cdot\overline{(\bE_h^g\times\bn)}) ds +\alpha(\nabla\sigma_h^k,\nabla d_k)_{\Omega_c}}{\|\bE_h^g\times\bn\|^2_{L^2(\Gamma)}+\alpha\|\nabla d_k\|^2_{L^2(\Omega_c)}};$$
\item[7.] Update 
$\sigma^{k+1}_h=\sigma_h^k+\gamma_k d_k$; set $k:=k+1$ and go to Step 1 until convergence is achieved.
\end{itemize}
\end{alg}
We note that the step size $\gamma_k$ in Step 6 above 
is not calculated by the exact line search algorithm, but it is simply computed 
by using the quadratic approximation of the objective function $\Phi_\alpha$ at $\sigma_0+\sigma_h^k$ along direction $d_k$, namely, the real quadratic function $\Psi(\gamma)$ in \eqref{linearap} with $\sigma_a=\sigma_h^k, \sigma_b=d_k$.  

\subsection{Sobolev gradient}
We recall that we have defined and used the weak gradient of the objective functional 
\eqref{obj-discrete} in \eqref{grad0d} or \eqref{grad_wd} that approximates the continuous 
gradient in \eqref{grad_w}.
It appears that the nonlinear conjugate Algorithm \ref{NLCG} converges very slowly 
for our nonlinear eddy current inverse problem, similarly to its behavior for most other nonlinear 
inverse problem. Next, we introduce a Sobolev gradient to help improve the convergence 
as it was done in \cite{JZ}. We can easily see that  the weak gradient $g(\sigma)$ in \eqref{grad_w} is just 
the weak gradient of the objective functional in the $L^2$ sense. Now we define 
a Sobolev gradient of the functional in the $H^1$ sense, namely to find an element 
$g_S(\sigma)\in H^1_0(\Omega_c)$ satisfying 
\begin{eqnarray}
\int_{\Omega_c}\nabla g_S(\sigma)\cdot\nabla\psi + g_S(\sigma)\psi dx = \int_{\Omega_c}g(\sigma)\psi dx ~~\forall\, \psi \in H^1_0(\Omega_c),\label{sobog}
\end{eqnarray}  
which is the weak formulation of the elliptic eproblem   
 \begin{eqnarray*}
\left\{
\begin{array}{rlll}
-\Delta g_S(\sigma) + g_S(\sigma) &=& g(\sigma)&\mbox{ in } \Omega_c,\\
g_S(\sigma) &=& 0& \mbox{ on } \partial\Omega_c.
\end{array}\right.
\end{eqnarray*}
This suggests us to compute the Sobolev gradient $g^S_h$ in the third step of Algorithm \ref{NLCG} 
by solving the following equation
\begin{eqnarray}
&&\int_{\Omega_c}\nabla g^S_h \cdot \nabla\tilde\sigma_h +g^S_h \tilde\sigma_h dx \nonumber\\
&=&
\alpha\int_{\Omega_c}\nabla\sigma_h\cdot\nabla\tilde\sigma_h dx 
-\omega\int_{\Omega_c}(\bE_h^1\cdot\bF_h^2-\bE_h^2\cdot\bF_h^1)\tilde\sigma_h dx ~~\forall\,\tilde\sigma_h\in V_h .\label{sobograd}
\end{eqnarray}
This can be solved very efficiently by many existing preconditioning-type iterative methods. 
\section{Numerical experiments}\label{sect5}
In this section, we present some numerical examples to illustrate the efficiency of Algorithm \ref{NLCG}. We take 
the computational domain $\Omega=[-2,2]\times[-2,2]\times[-2,0.2]$, with the non-conducting subregion  
$\Omega_0=[-2,2]\times[-2,2]\times[0,0.2]$ (where the conductivity $\sigma_0$ vanishes) and 
the conducting subregion $\Omega_c=[-2,2]\times[-2,2]\times[-2,0]$.
The state and adjoint state equations involved are solved with edge element methods. 
We implement the algorithm using the parallel hierarchical grid platform (PHG) \cite{phg}. 
The numerical examples are carried out using an Apple laptop with Intel i7 8750h  CPU and 16G memory.  
The data $\bn\times\bE^{obs}$ is generated by the edge element method \cite{CCCZ}, and 
can be written as 
$$\bn\times\bE^{obs}(x)|_\Gamma=\bn\times\sum_{e\in\mathcal{E}_\Gamma} (R^e + i I^e)  \Phi^e(x),$$ 
where $\mathcal{E}_\Gamma$ is the union of all edges of the mesh $\mathcal{M}_h$ on the measurement surface 
$\Gamma$ (i.e., the plane $z=0.2$), 
$R^e$ and $I^e$ are the degrees of freedom on the edge $e$ for real and imaginary parts of the electric field with exact abnormal conductivity, respectively, and $\Phi^e(x)$ is the edge element basis function 
associated with edge $e$. 
To test the algorithm with noisy data,  we generate the noisy data by adding the noise in the form
\begin{eqnarray*}
\bn\times\bE^{obs}(x)|_\Gamma=\bn\times\sum_{e\in\mathcal{E}_\Gamma} (R^e + i I^e)(1+\delta\xi)\Phi^e(x),
\end{eqnarray*}
where $\delta$ is the noise level,  and $\xi$ is a uniformly distributed random variable in $[-1,1]$. 

In all examples we choose the source $\bJ_s$ as ($\nabla\cdot\bJ_s=0$) 
\begin{eqnarray*}
\bJ_s = \nabla\times \sum^9_{i,j=1}\delta(\bx-\bx_{ij})\be_1,
\end{eqnarray*}
where $\be_1$ is the unit vector alone the $x$-axis and 
$\bx_{ij}=(-2.0+0.4*i,-2.0+0.4*j,0.1)$, i.e., there are 81 point sources on plane $z=0.1$. We assume the background conductivity $\sigma_0=1.0$ in $\Omega_c=\Omega_1\cup\Omega_2$.
By this setting, we apply Algorithm \ref{NLCG} to recover the abnormal conductivity $\sigma$ with the data on boundary $\Gamma$. 
We always choose the initial guess $0$ in the NLCG algorithm, and take the parameters 
$\varepsilon=1.0,\mu=1.0,\omega=0.79$, and the regularization $\alpha=10^{-6}$ unless 
it is specified otherwise.

\subsection{Example 1}
In this example, the domain with abnormal conductivity is 
$\Omega_2=[-0.4,0.4]\times[-0.4,0.4]\times[-1.2,-0.4]$, 
where the exact abnormal conductivity is given by $\sigma=1.0$, and $\sigma$ vanishes in $\Omega_1$.  
That is, the exact conductivity
 $\sigma_0+\sigma$ is constant $0.0$, $1.0$ and $2.0$ in $\Omega_0$, $\Omega_1$ and $\Omega_2$ 
 respectively. 

\begin{figure}
\begin{center}
\includegraphics[width=0.45\textwidth]{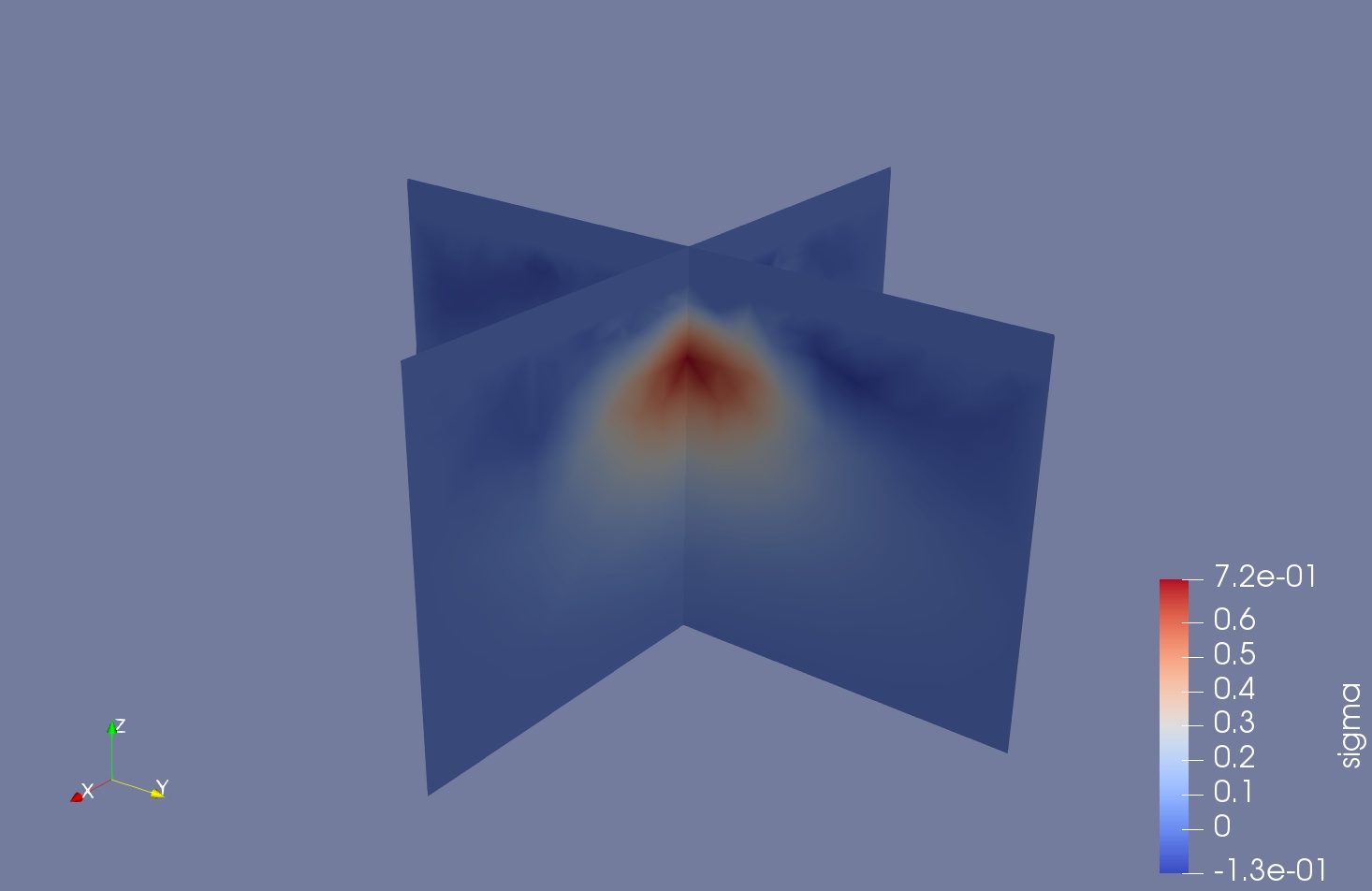}
\includegraphics[width=0.45\textwidth]{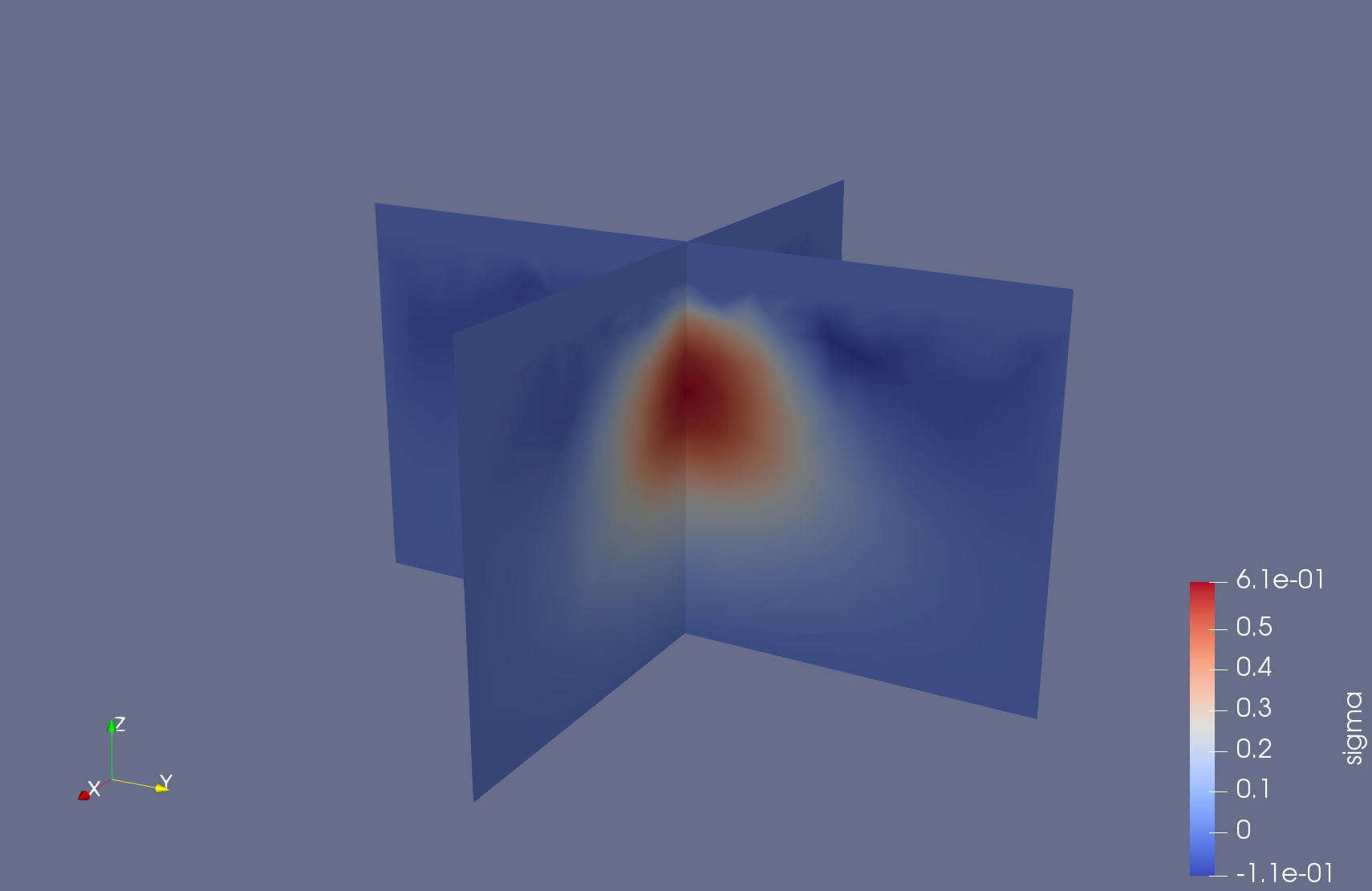}
\includegraphics[width=0.45\textwidth]{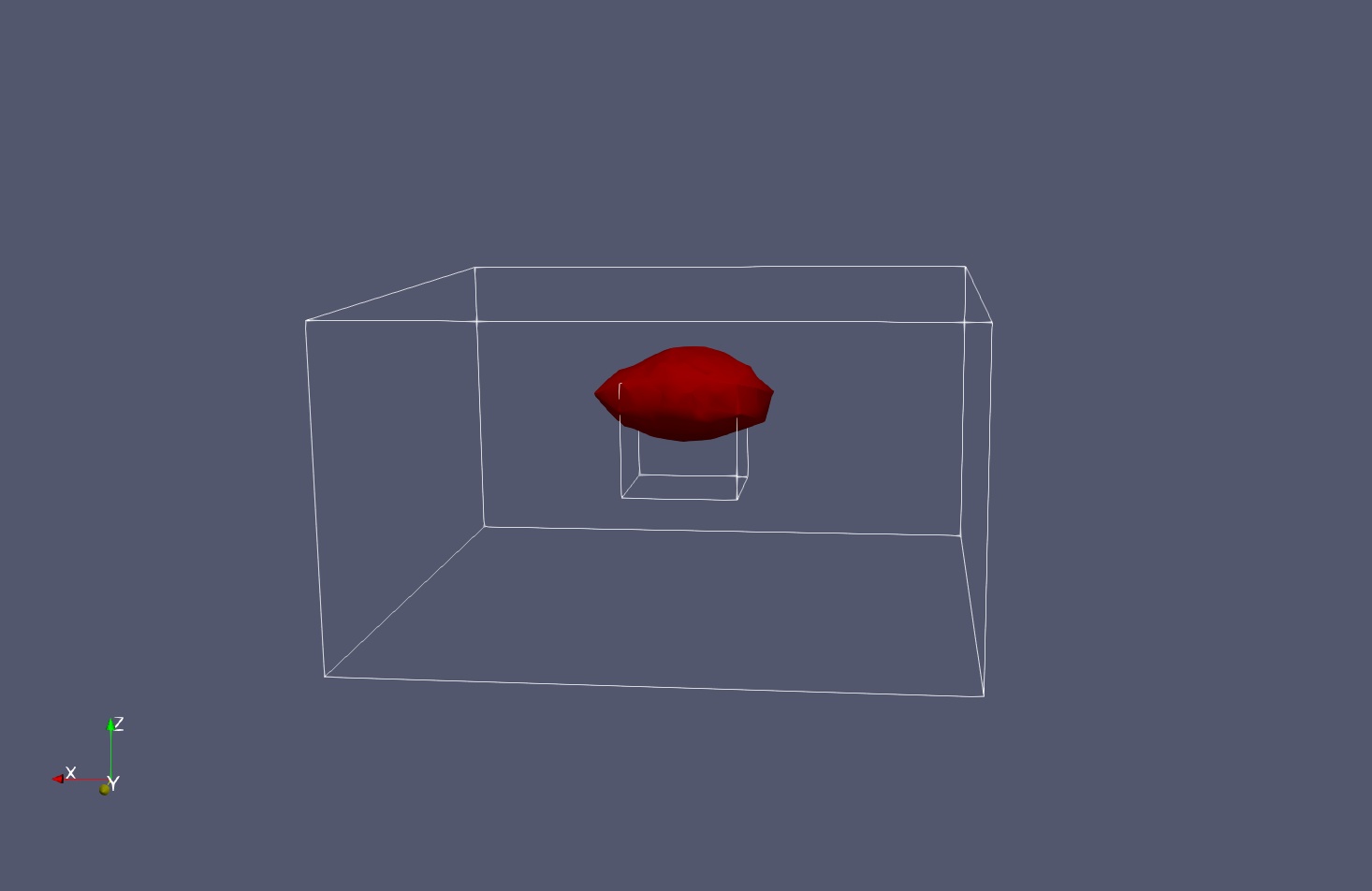}
\includegraphics[width=0.45\textwidth]{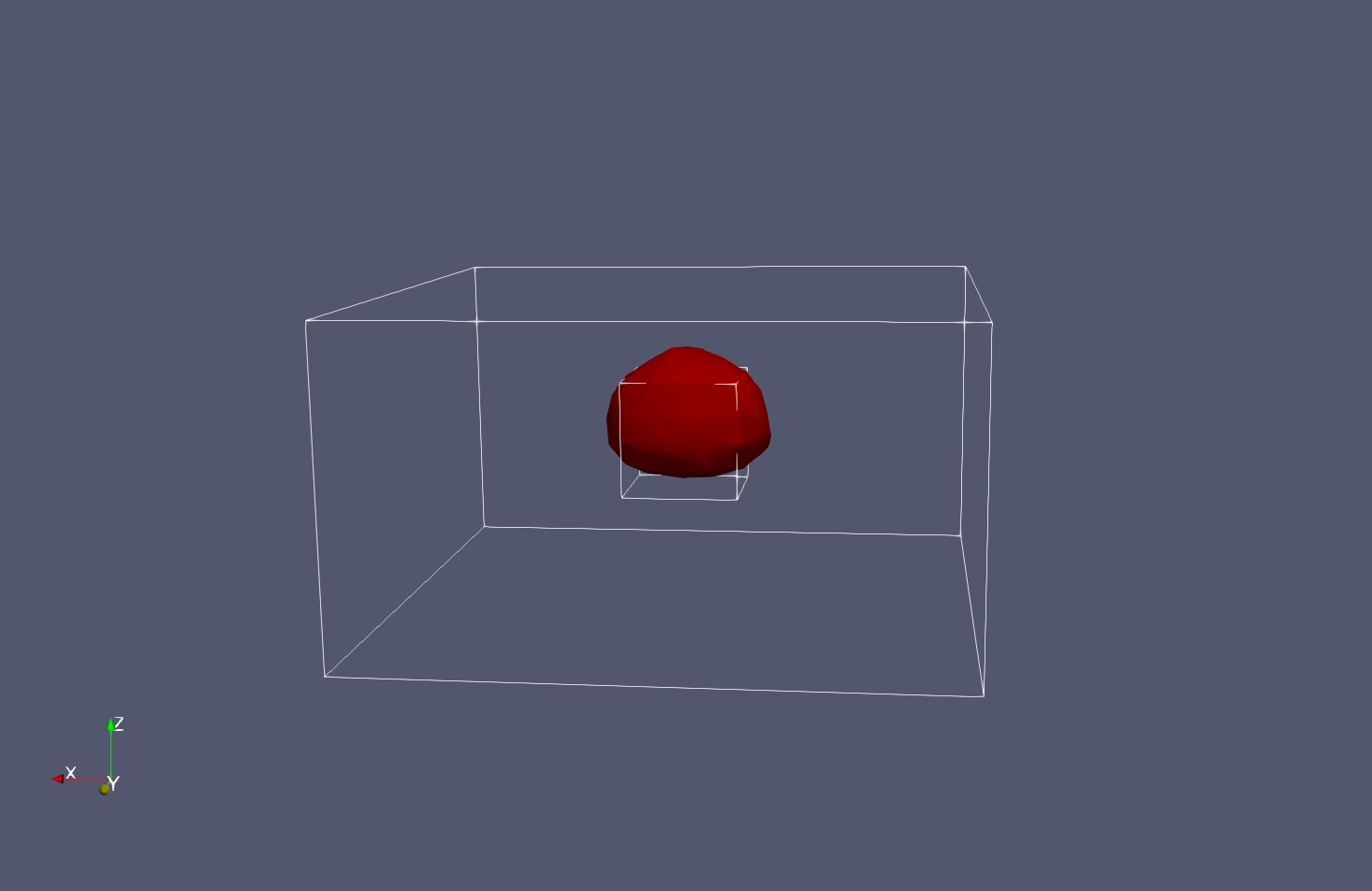}
\caption{The recovery of $\sigma$ after 100 iterations (Upper left) and 200 iterations (Upper right); The lower two pictures are the corresponding isosurfaces of the recovered $\sigma$ with isovalue 0.35.}\label{fig2}
\end{center}
\end{figure}

\begin{figure}
\begin{center}
\includegraphics[width=0.45\textwidth]{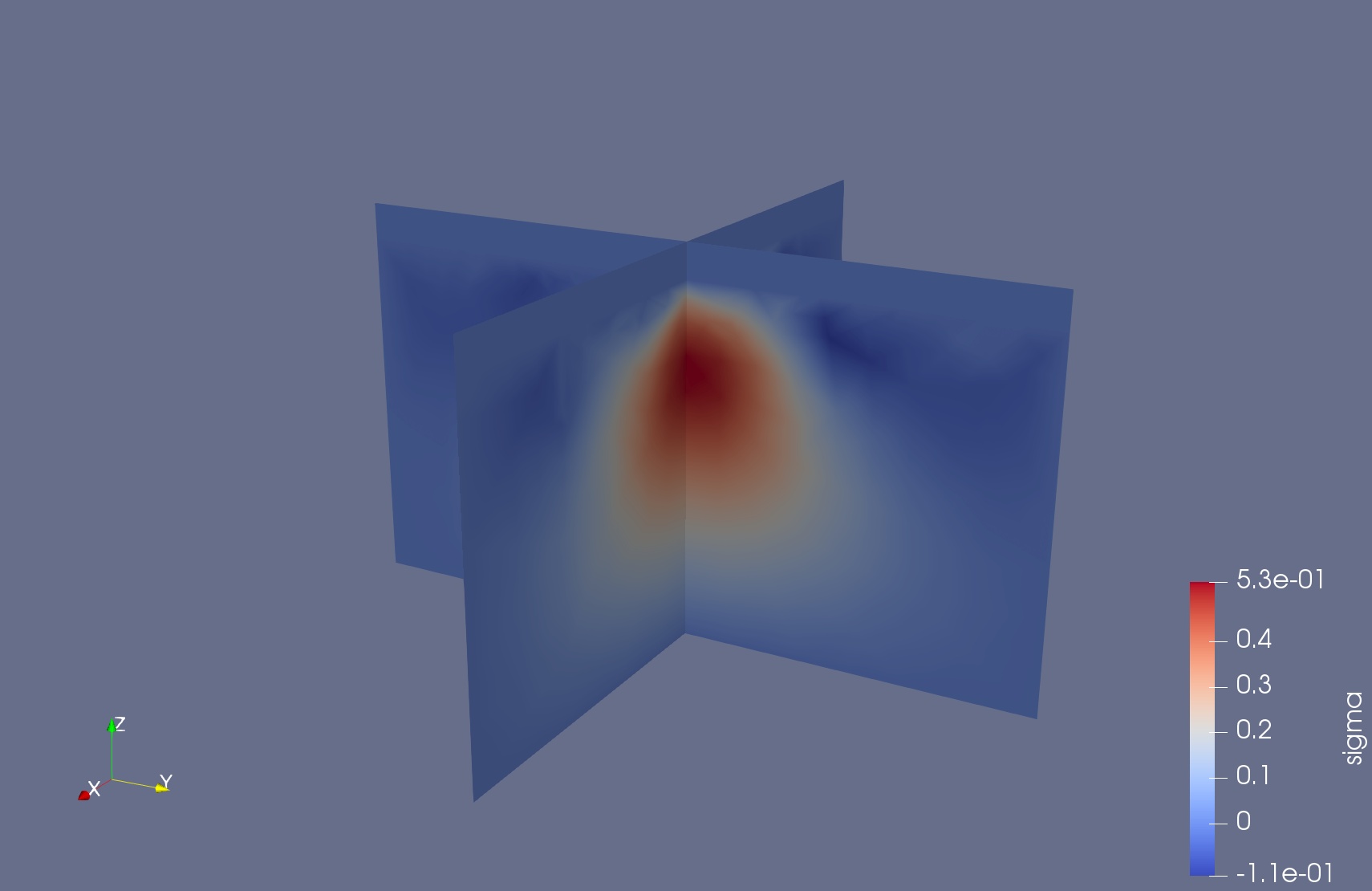}
\includegraphics[width=0.45\textwidth]{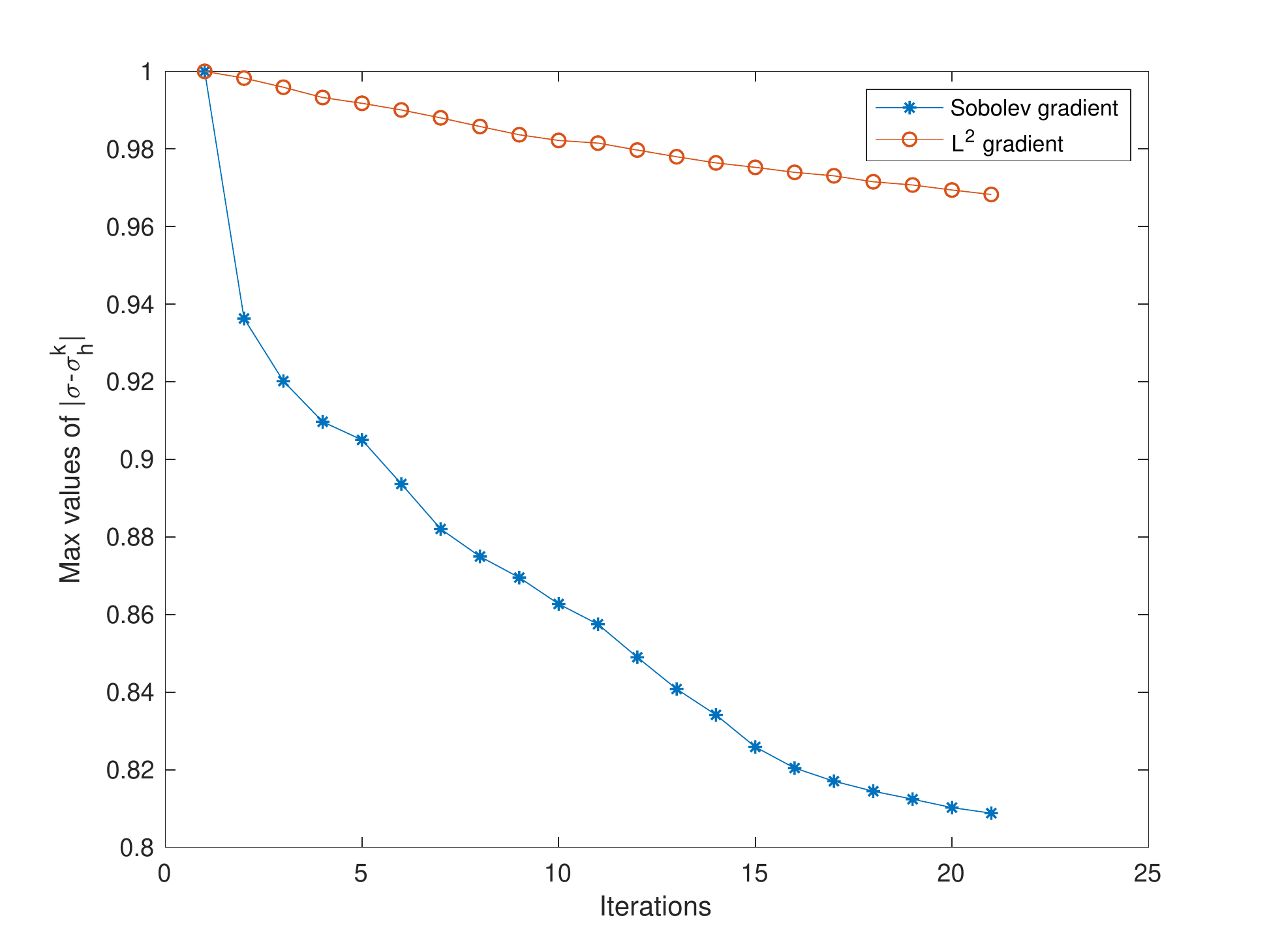}
\caption{The recovery of $\sigma$ after 20 iterations (Left) and the convergence history (Right).}\label{fig3}
\end{center}
\end{figure}

The total degrees of freedom of the edge elements are 213,128.  First, we use the $L^2$ gradient of the objective functional in Algorithm \ref{NLCG}, and the recovery results are shown in 
Figure \ref{fig2}, where the left and right pictures give the results in 100 and 200 iterations, 
respectively. 
We can find that the recovery is closer to the exact conductivity with more iterations.
The recovery result by Algorithm \ref{NLCG} using the Sobolev gradient is given in the left of 
Figure
\ref{fig3}
(20 iterations), with the convergence history of the nonlinear CG algorithm by 
the $L^2$ and Sobolev gradients, respectively, in the right of Figure \ref{fig3}. 
We notice that the algorithm with the Sobolev gradient converges much faster.  
The recovered conductivity by using Sobolev gradient with 20 iterations is very close 
to the result with 200 iterations by using the $L^2$ gradient.

\subsection{Example 2}
\begin{figure}
\begin{center}
\includegraphics[width=0.4\textwidth]{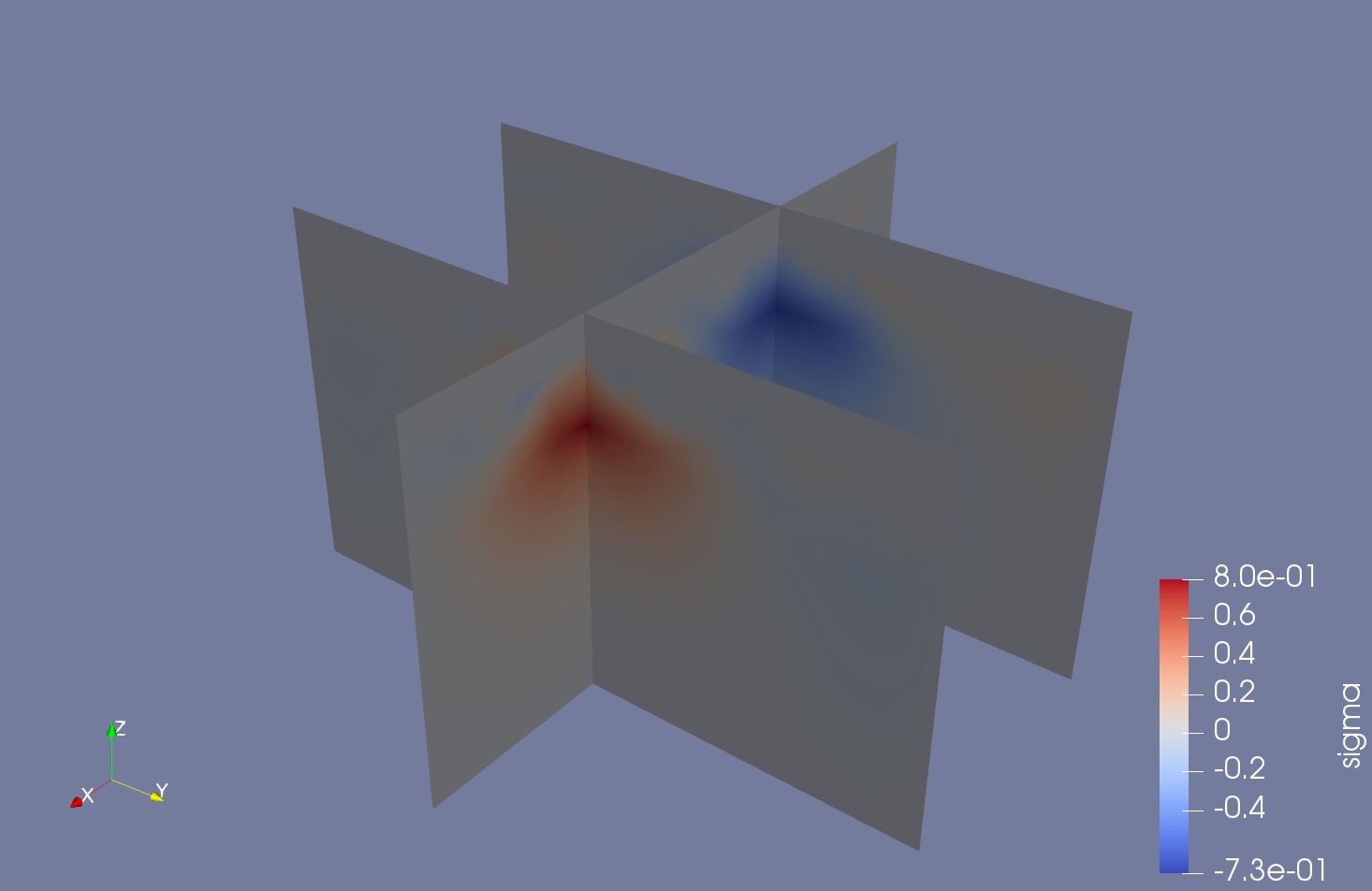}
\includegraphics[width=0.4\textwidth]{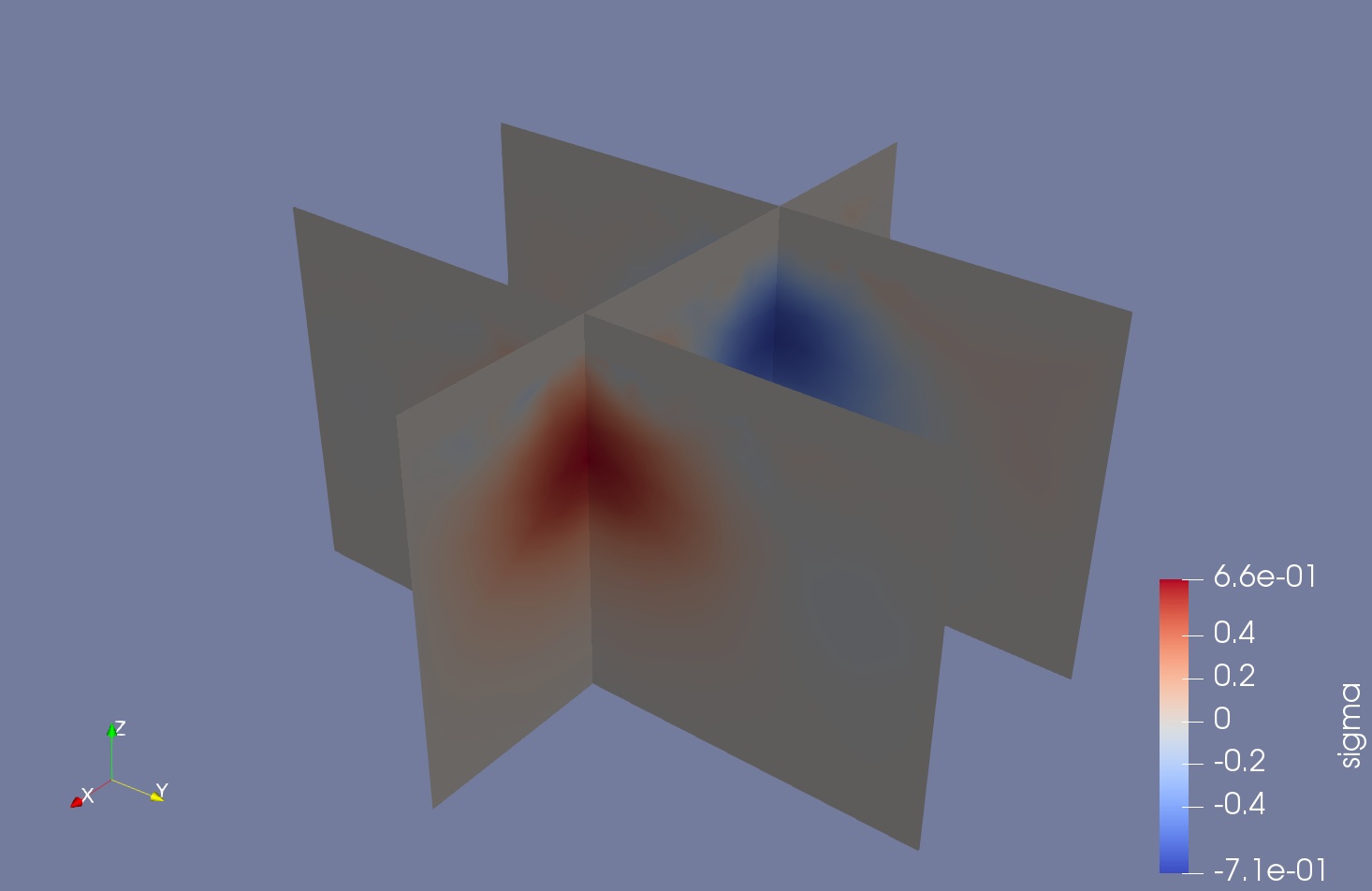}
\includegraphics[width=0.4\textwidth]{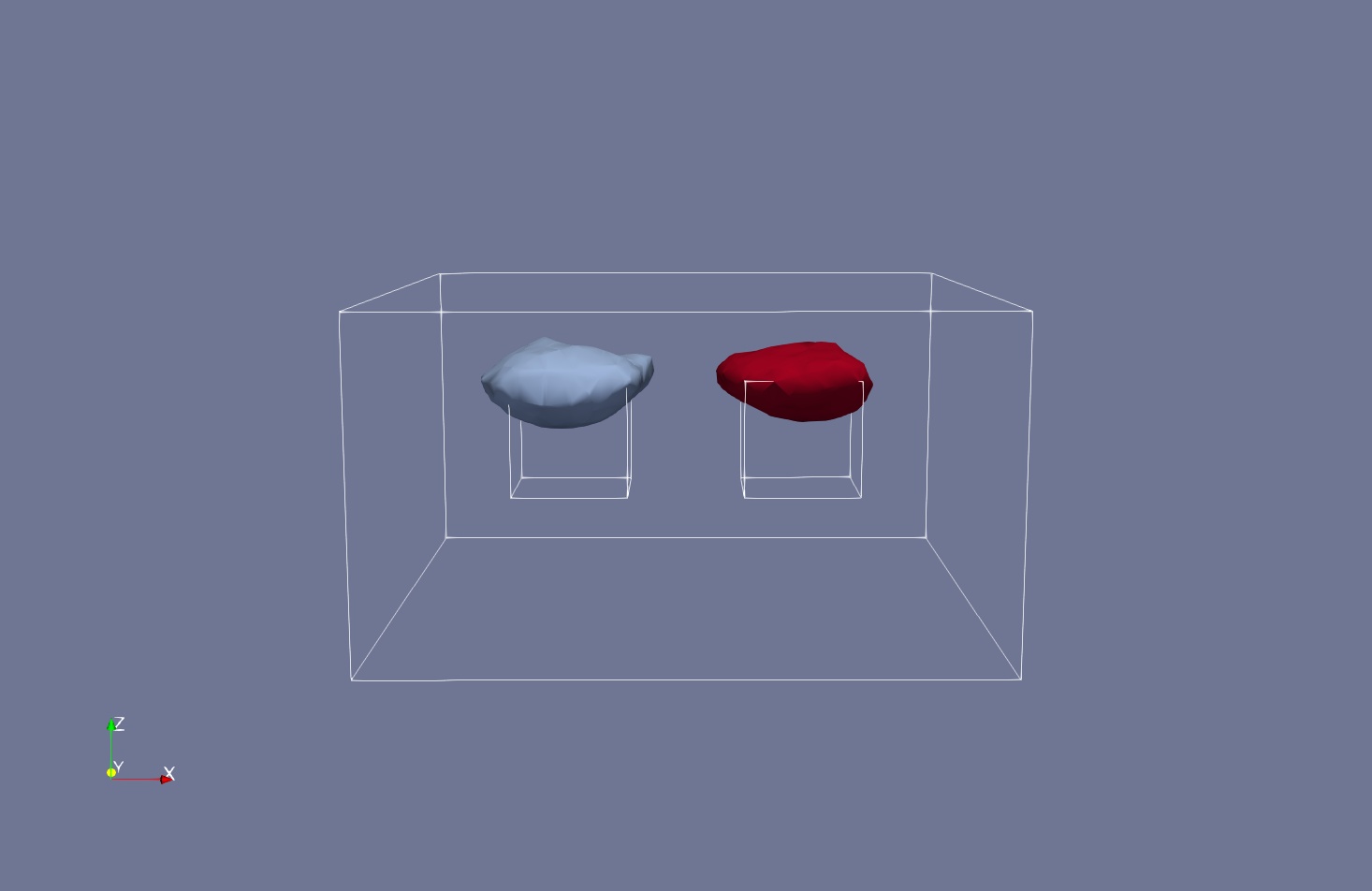}
\includegraphics[width=0.4\textwidth]{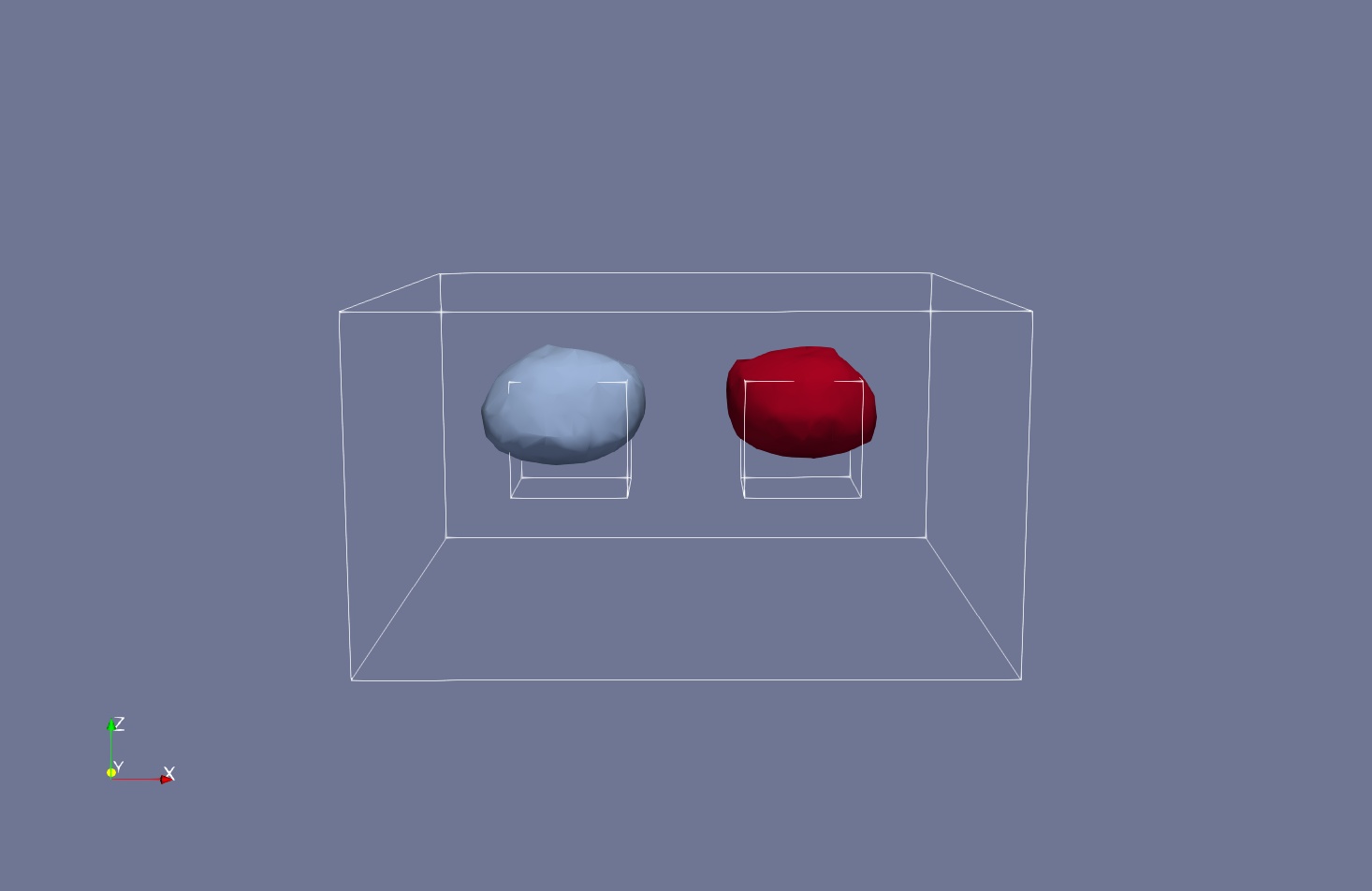}
\caption{The recovered $\sigma$ after 100 iterations (Upper left) and 200 iterations (Upper right). The lower two pictures are the isosufaces of the recovered $\sigma$ with isovalues 0.35 (Lower left) and -0.35 (Lower right). 
The small cubes are the real locations of the two anomalies.}\label{fig5}
\end{center}
\end{figure}

 In this example, we consider the case with two abnormal subdomains, $\Omega_{21}$ and $\Omega_{22}$. 
 We take $\Omega_{21}=[-1.2,-0.4]\times[-0.4,0.4]\times[-1.2,-0.4]$, with the exact abnormal conductivity $\sigma=-0.9$, and $\Omega_{22}=[0.4,1.2]\times[-0.4,0.4]\times[-1.2,-0.4]$, with the exact abnormal conductivity $\sigma=1.0$.
To be more specific,
\begin{eqnarray*}
\sigma_0+\sigma=\left\{\begin{array}{cc} 0 &\mbox{ in } \Omega_0,\\
                             1.0 &\mbox{ in } \Omega_1,\\
                             0.1 &\mbox{ in } \Omega_{21},\\
                             2.0 &\mbox{ in } \Omega_{22}.
                             \end{array}\right.
\end{eqnarray*}
The total degrees of freedom of the edge elements of the state and adjoint equations are 266,690. 
In this example, we show the recovery results only for the Sobolev gradient defined by \eqref{sobograd}. 
The left and right pictures of Figure\,\ref{fig5} present the recovery results in 100 and 200 
iterations, respectively. We can see that two abnormal objects are well separated and recovered. 

\begin{figure}
\begin{center}
\includegraphics[width=0.4\textwidth]{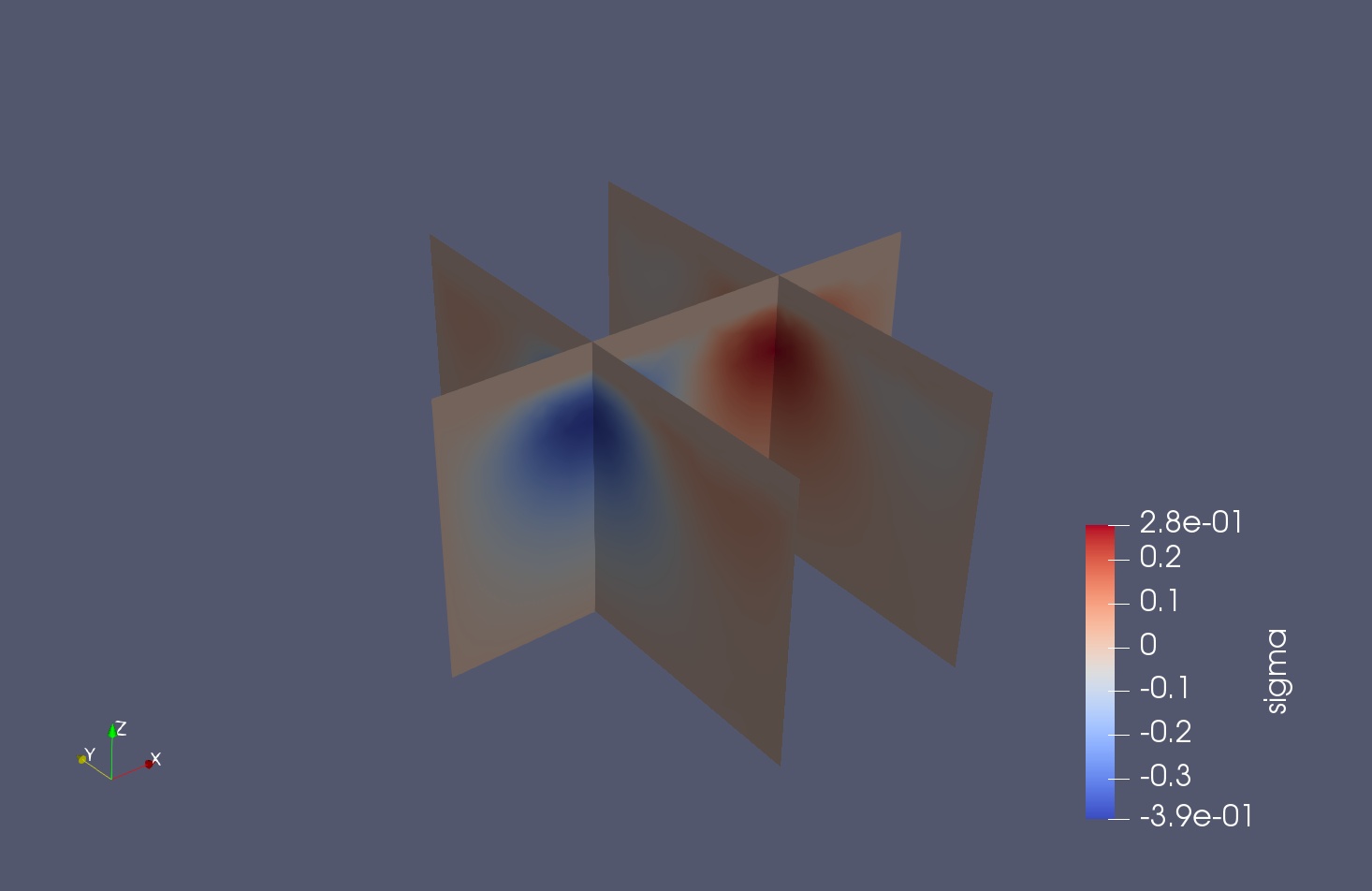}
\includegraphics[width=0.4\textwidth]{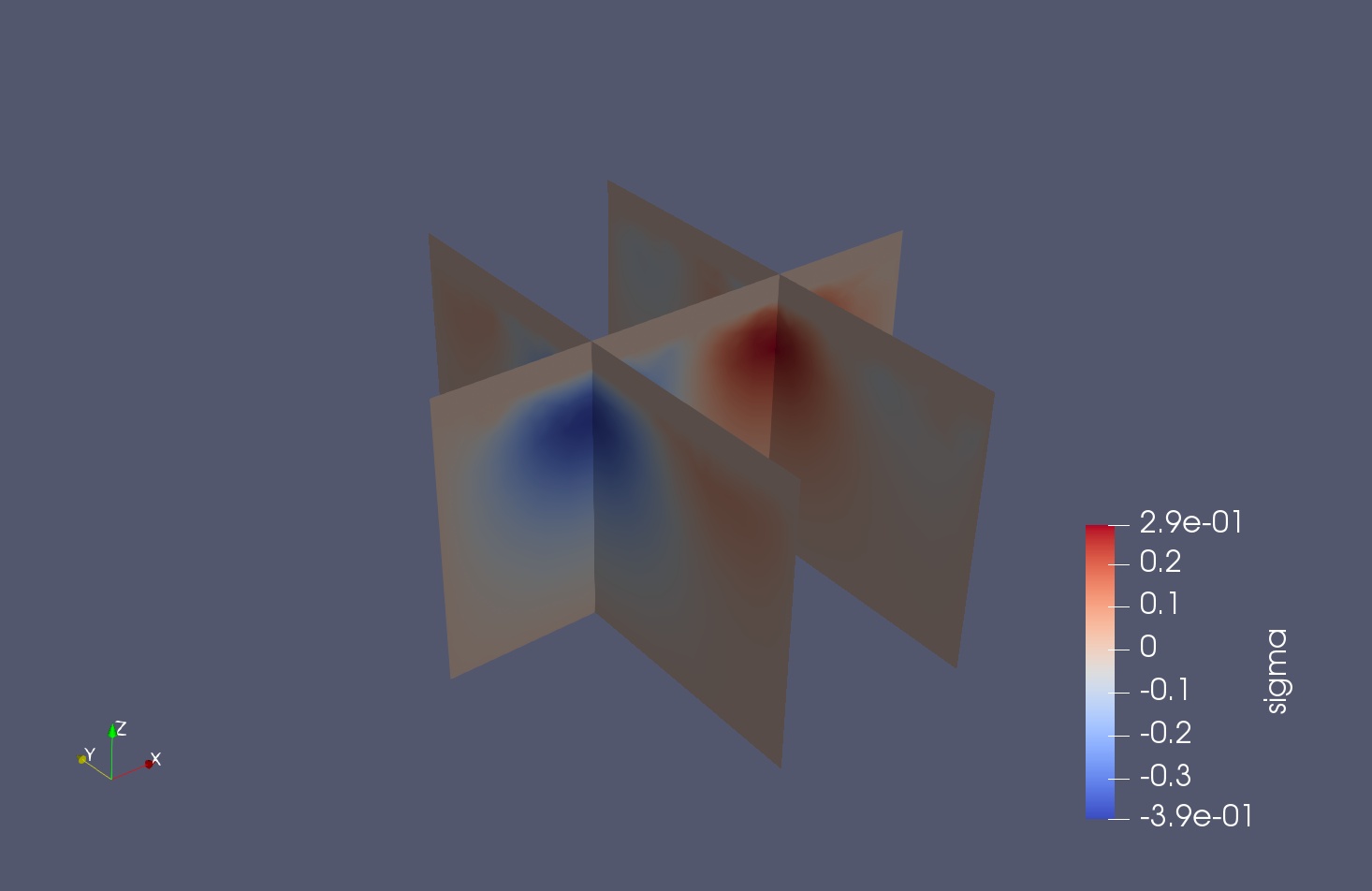}
\includegraphics[width=0.4\textwidth]{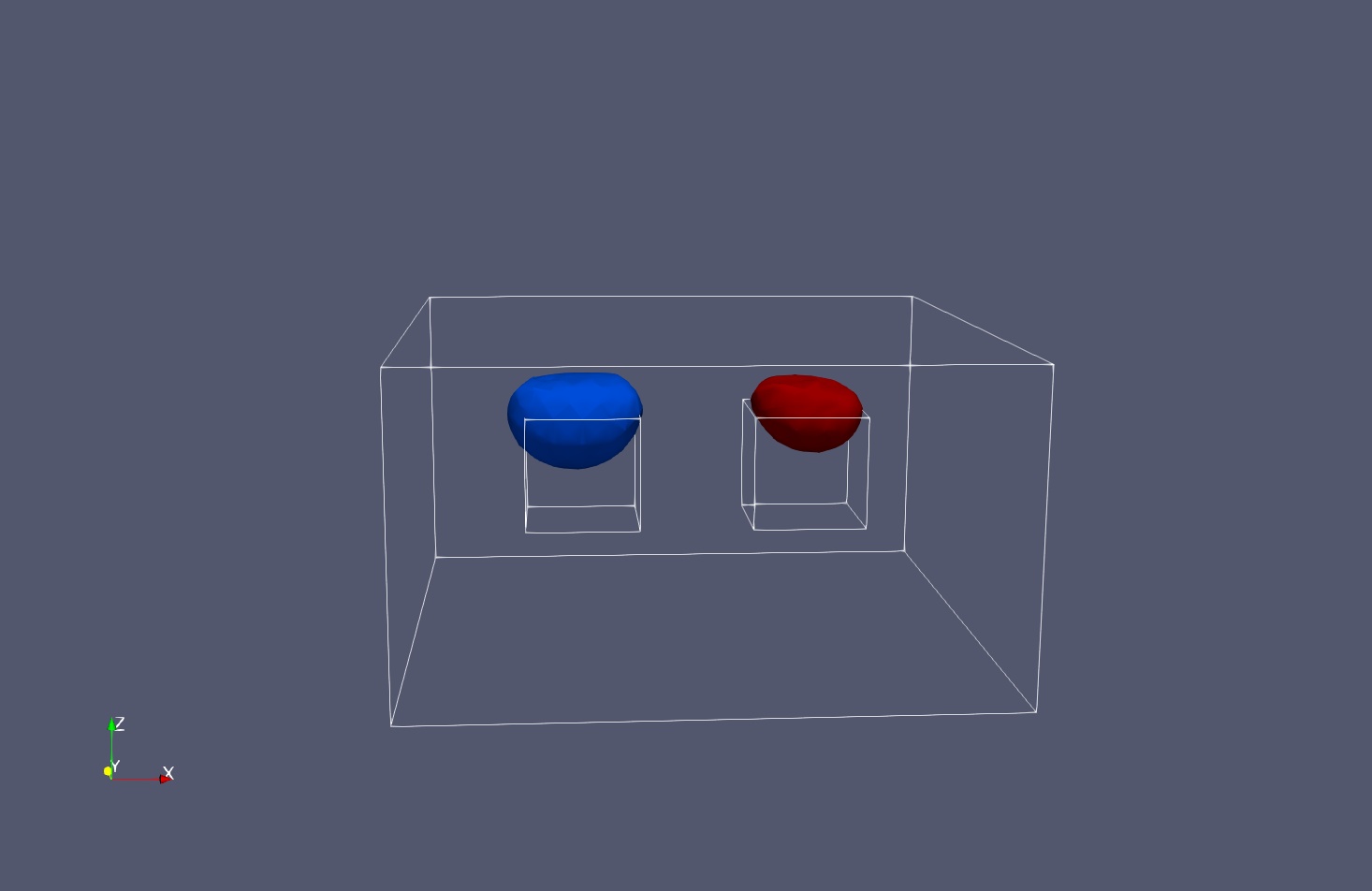}
\includegraphics[width=0.4\textwidth]{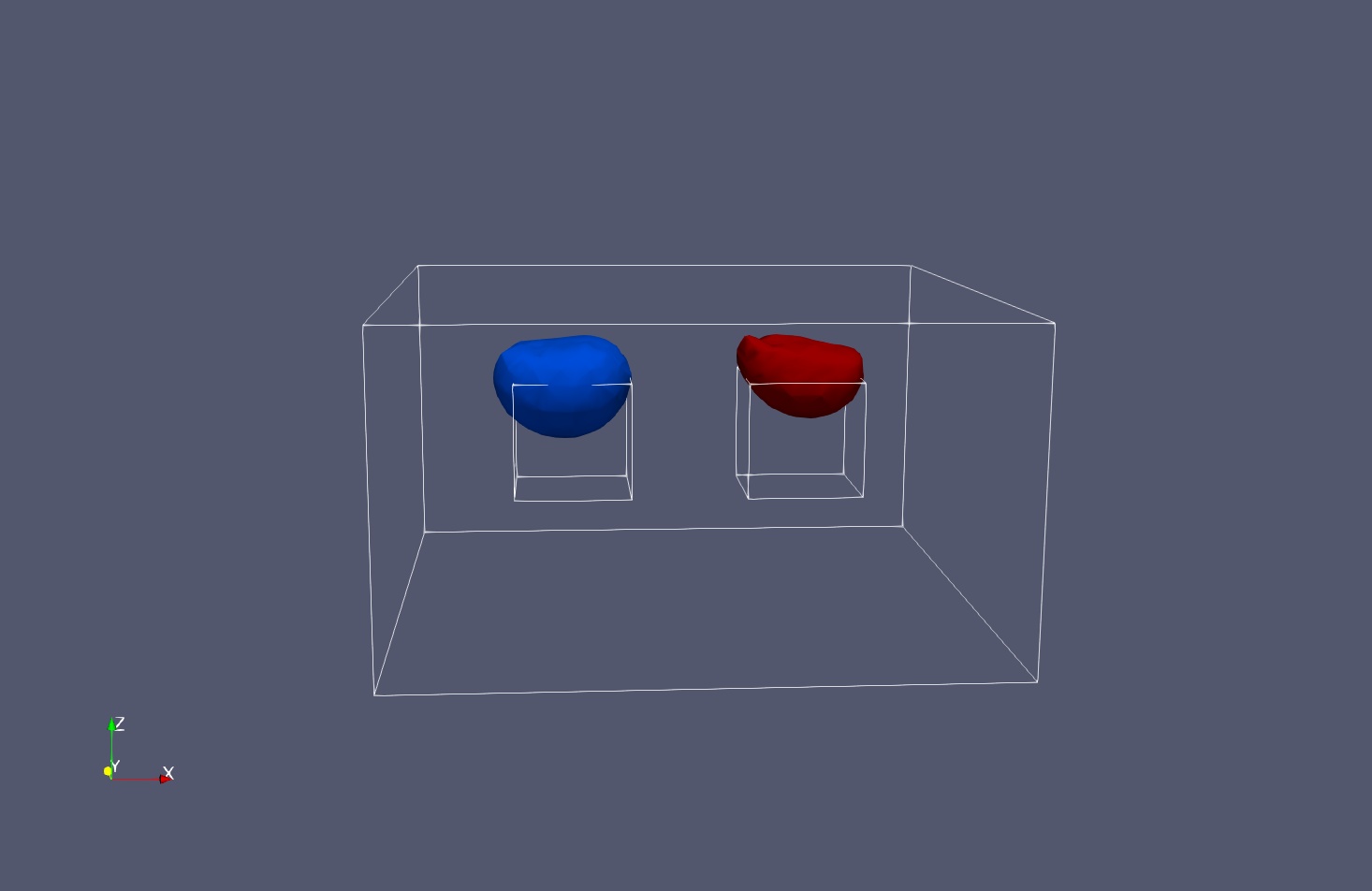}
\caption{Recovery results on after 100 iterations with noisy data, the noise level is $0.1\%$ (Upper left)  and  $0.4\%$ (Upper right).
The lower two pictures are the isosufaces of recovered $\sigma$ with isovalues 0.20 (Lower left) and -0.25 (Lower right).}\label{fig6}
\end{center}
\end{figure}

Figure \ref{fig6} shows the recovered conductivity with the noisy data: 
the left picture with the noise level $\delta=0.1\%$, and the right picture with $\delta=0.4\%$, for both of which 
the regularization parameter is taken to be $\alpha = 10^{-4}$ and 100 iterations are conducted. 
\section{Concluding remarks}\label{sect6}

We have studied an ill-posed eddy current inversion problem mathematically and numerically. 
We have first investigated the ill-posedness of the inverse eddy current problem, by showing 
the compactness of the forward operator mapping the conductivity parameter to the tangential trace of the electric field, and the non-uniqueness of the inverse problem. 
For the nonlinear regularized minimization formulation of the inverse problem, we have explored 
the existence and stability of the minimizers, and the optimality system of its Lagrange formulation
in terms of the real and imaginary parts of the PDE constraints, as well as the finite element 
approximation of the nonlinear regularized minimization system. A nonlinear conjugate gradient method 
is formulated for solving the discrete nonlinear constrained optimization problem, with 
its step sizes updated very effectively by a quadratic approximation of the objective function, and a Sobolev gradient introduced to effectively accelerate the iteration. 
Numerical examples have shown the feasibility and effectiveness of the reconstruction algorithm, 
which can clearly recover the locations and sizes of separated inclusions in the noisy case. 

%
${\ }$

\end{document}